\documentclass[journal]{IEEEtran}

\usepackage{cite}
\usepackage{amsmath}
\usepackage{amssymb}
\usepackage{graphicx}
\usepackage{subfigure}
\usepackage{epstopdf}
\usepackage{amsfonts}
\usepackage[all]{xy}
\usepackage{array,xcolor}
\usepackage{color}

\usepackage{caption}
\usepackage{algorithm}
\usepackage{algpseudocode}
\usepackage{multirow}

\usepackage{threeparttable}
\usepackage{booktabs}
\usepackage{changepage}

\usepackage{stfloats}
\usepackage{balance}

\newtheorem{myDef}{Definition}
\newtheorem{myTheo}{Theorem}

\newtheorem{myRem}{Remark}
\newtheorem{myLem}{Lemma}
\newtheorem{myCor}{Corollary}
\newtheorem{myAss}{Assumption}

\begin{document}
\title{Platooning of Connected Vehicles with Undirected Topologies: Robustness Analysis and Distributed H-infinity Controller Synthesis}
\author{Yang~Zheng,~\IEEEmembership{Student Member,~IEEE}, Shengbo~Eben~Li,~\IEEEmembership{Senior Member,~IEEE}, Keqiang~Li, \\ and~Wei~Ren,~\IEEEmembership{Fellow,~IEEE}
\thanks{Y. Zheng is supported by the Clarendon Scholarship and the Jason Hu Scholarship. This work is partially supported by NSF China with grant 51622504 and 51575293, and National Key R\&D Program of China with 2016YFB0100906. \emph{(Corresponding author: Yang Zheng)}}
\thanks{Y. Zheng is with the Department of Engineering Science, University of Oxford, Parks Road, Oxford OX1 3PJ, U.K. (e-mail: yang.zheng@eng.ox.ac.uk.)}
\thanks{S. Li and K. Li are with the Department of Automotive Engineering, Tsinghua University, Beijing, 100084, China (e-mail: lisb04@gmail.com; likq@tsinghua.edu.cn)}
\thanks{W. Ren is is with the Department of Electrical and Computer Engineering, University of California, Riverside, CA 92521 USA. (e-mail: ren@ee.ucr.edu)}
}

\maketitle

\begin{abstract}
    This paper considers the robustness analysis and distributed $\mathcal{H}_{\infty}$ (H-infinity) controller synthesis for a platoon of connected vehicles with undirected topologies. We first formulate a unified model to describe the collective behavior of homogeneous platoons with external disturbances using graph theory. By exploiting the spectral decomposition of a symmetric matrix, the collective dynamics of a platoon is equivalently decomposed into a set of subsystems sharing the same size with one single vehicle. Then, we provide an explicit scaling trend of robustness measure $\gamma$-gain, and introduce a scalable multi-step procedure to synthesize a distributed $\mathcal{H}_{\infty}$ controller for large-scale platoons. It is shown that communication topology, especially the leader's information, exerts great influence on both robustness performance and controller synthesis. Further, an intuitive optimization problem is formulated to optimize an undirected topology for a platoon system, and the upper and lower bounds of the objective are explicitly analyzed, which hints us that coordination of multiple mini-platoons is one reasonable architecture to control large-scale platoons. Numerical simulations are conducted to illustrate our findings.
\end{abstract}

\begin{IEEEkeywords}
Connected vehicles, platoon control, robustness analysis, distributed $\mathcal{H}_{\infty}$ control, topology design.
\end{IEEEkeywords}

\IEEEpeerreviewmaketitle

\section{Introduction}

\IEEEPARstart{T}he increasing traffic demand in today's life brings a heavy burden on the existing transportation infrastructure, which sometimes leads to a heavily congested road network and even results in serious casualties. Human-centric methods to these problems provide some insightful solutions, but most of them are constrained by human factors, \emph{e.g.}, reaction time and perception limitations~\cite{dong2011driver,wang2015driving,zheng2014driving}. On the other hand, vehicle automation and multi-vehicle cooperation are very promising to enhance traffic safety, improve traffic capacity, and reduce fuel consumption, which attracts increasing attention in recent years (see \cite{naus2010string,oncu2014cooperative,zheng2015dynamic} and the references therein).

The platooning of connected vehicles, an importation application of multi-vehicle cooperation, is to ensure that all the vehicles in a group maintain a desired speed and keep a pre-specified inter-vehicle spacing. The platooning practices date back to the PATH program during the last eighties \cite{shladover1991automated}. Since then, many topics on platoon control have been addressed, such as the selection of spacing policies \cite{swaroop1994comparision}, the influence of imperfect communication \cite{liu2001effects, fernandes2012platooning}, and the impacts of heterogeneity on string stability \cite{shaw2007string,lestas2007scalability}. Recently, advanced control methods have been introduced and implemented in order to achieve better performance for platoons: Dunbar and Derek (2012) introduced a distributed receding horizon controller for platoons with predecessor-following topology \cite{dunbar2012distributed}, which is recently extended to unidirectional topologies in~\cite{Zheng2016distributed}; Ploeg \emph{et al.} (2014) proposed an $\mathcal{H}_{\infty}$  controller synthesis approach for platoons with linear dynamics, where string stability was explicitly satisfied by solving a linear inequality matrix (LMI) \cite{ploeg2014controller}; Zheng \emph{et al.} (2016) explicitly derived the stabilizing thresholds of the controller gains by using the graph theory and Routh-Hurwitz criterion, which could cover a large class of communication topologies \cite{zheng2016stability}; Zhang and Orosz (2016) introduced a motif-based approach to investigate the effects of heterogeneous connectivity structures and information delays on platoon systems~\cite{zhang2016motif}. The interested reader can refer to a recent review in \cite{li2015overview}.

One recent research focus is on finding essential performance limitations of large-scale platoons~\cite{zheng2015dynamic, li2015overview,seiler2004disturbance,barooah2005error,middleton2010string,barooah2009mistuning, hao2011stability,zheng2015stabilityMargin}. Many works focused on two kinds of performance measures: 1) string stability, which refers to the attenuation effect of spacing error along the vehicle string \cite{seiler2004disturbance}; 2) stability margin, which characterizes the convergence speed of initial errors~\cite{zheng2015stabilityMargin}. For example, Seiler \emph{et al.} (2004) proved that string stability cannot be satisfied for homogeneous platoons with predecessor-following topology and constant spacing policy due to a complementary sensitivity integral constraint~\cite{seiler2004disturbance}. Barooah \emph{et al.} (2005) showed that platoons with bidirectional topology also suffered fundamental limitations on string stability~\cite{barooah2005error}. Middleton and Braslavsky (2010) further pointed out that both forward communication and small time-headway cannot alter the limitations on string stability for platoons, in which heterogeneous vehicle dynamics and limited communication range were considered~\cite{middleton2010string}. As for the limitations of stability margin in a large-scale platoon, Barooah \emph{et al.} (2009) proved that the stability margin would approach zero as $\mathcal{O}(1/N^2 )$  ($N$ is the platoon size) for symmetric control, and demonstrated that the asymptotic behavior could be improved to $\mathcal{O}(1/N)$  via introducing certain "mistuning"~\cite{barooah2009mistuning}. Using partial differential equation (PDE) approximation, Hao \emph{et al.} (2011) showed that the scaling of stability margin could be improved to $\mathcal{O}(1/N^{2/D} )$ under D-dimensional communication topologies \cite{hao2011stability}. Zheng \emph{et al.} (2016) further introduced two useful methods, \emph{i.e.}, enlarging the communication topology and using asymmetric control, to improve the stability margin from the perspective of topology selection and control adjustment in a unified framework \cite{zheng2015stabilityMargin}. These studies have offered insightful viewpoints on the performance limitations of large-scale platoons in terms of string stability and stability margin.

In this paper, we focus on the robustness analysis and controller synthesis of large-scale platoons with undirected topologies considering external disturbances. This paper shows additional benefits on understanding the essential limitations of platoons, and also provides a distributed $\mathcal{H}_{\infty}$ method to design the controller with guaranteed performance. Using algebraic graph theory, we first derive a unified model in both time and frequency domain to describe the collective behavior of homogeneous platoons with external disturbances. A $\gamma$-gain is used to quantify the robustness of a platoon from the perspective of energy amplification. The major strategy of this paper is to equivalently decouple the collective dynamics of a platoon into a set of subsystems by exploiting the spectral decomposition of a symmetric matrix \cite{ren2008distributed,li2011h}. Along with this idea, both robustness analysis and controller synthesis are carried out based on the decomposed subsystems which share the same dimension with one single vehicle. This fact not only significantly reduces the complexity in analysis and synthesis, but also explicitly highlights the influence of communication topology and the importance of leader's information. The contributions of this paper are:

\begin{enumerate}
  \item We analytically provide the scaling trend of robustness measure $\gamma$-gain for platoons with undirected topologies, which is lower bounded by the minimal eigenvalue of a matrix associated with the communication topology. Besides, we prove that $\gamma$-gain increases at least as $\mathcal{O}(N)$, if the number of followers that are pinned to the leader is fixed. For platoons with bidirectional topology, the scaling trend is deteriorated to $\mathcal{O}(N^2)$. These results provide new understandings on the essential limitations of large-scale platoons, which are also consistent with previous results on stability margin~\cite{barooah2009mistuning, zheng2015stabilityMargin, zheng2016stability}.
  \item We introduce a scalable multi-step procedure to synthesize a distributed $\mathcal{H}_{\infty}$ controller for large platoons. This problem can be equivalently converted into a set of $\mathcal{H}_{\infty}$ control of independent systems that share the same dimension with a single vehicle. It is shown that the existence of a distributed $\mathcal{H}_{\infty}$ controller is independent of the topology as long as there exists a spanning tree in the undirected topology.
  \item We give useful discussions on the selections of communication topologies based on the analytical results above. An intuitive optimization problem is formulated to optimize an undirected topology for a platoon system, where the upper and lower bounds of the objective are explicitly analyzed. These results hint us that a star topology might be a good communication topology considering limited communication resources, and that coordination of multiple mini-platoons is one reasonable architecture for the control of large-scale platoons.
\end{enumerate}

The rest of this paper is organized as follows. The problem statement is introduced in Section~\ref{section:pro}. Section~\ref{section:Robustness} presents the results on robustness analysis. In Section~\ref{section:controller}, we give the distributed $\mathcal{H}_{\infty}$ controller synthesis of platoons, and discuss the design of communication topology.  This is followed by numerical experiments in Section~\ref{section:simulation}, and we conclude the paper in Section~\ref{section:conclusion}.

\emph{Notations}: The real and complex numbers are denoted by $\mathbb{R}$ and $\mathbb{C}$, respectively. We use $\mathbb{R}^n$ to denote the $n$-dimensional Euclidean space. Let $\mathbb{R}^{m\times n}$ be the set of $m\times n$ real matrices, and $\mathbb{S}^n$ be the symmetric matrices of dimension $n$. $I$ represents the identity matrix with compatible dimensions. $\text{diag}\{a_1,\ldots,a_n\}$ is a block diagonal matrix with $a_i, i = 1, \dots, N,$ as the main diagonal entries. We use $\lambda_i(A)$ to represent the $i$-th eigenvalue of a symmetric matrix $A \in \mathbb{S}^{n}$, and denote $\lambda_{\min}(A), \lambda_{\max}(A)$ as its minimum and maximal eigenvalue, respectively. Given a symmetric matrix $ X \in \mathbb{S}^n $, $ X \succ (\prec)\text{ } 0 $ means that $ X $ is positive (negative) definite. Give $A \in \mathbb{R}^{m\times n}$ and $B \in \mathbb{R}^{p\times q}$, we use $A\otimes B$ to denote the Kronecker product of $A$ and $B$.

\begin{figure}[!t]
    \centering
    \includegraphics[width=0.95\columnwidth]{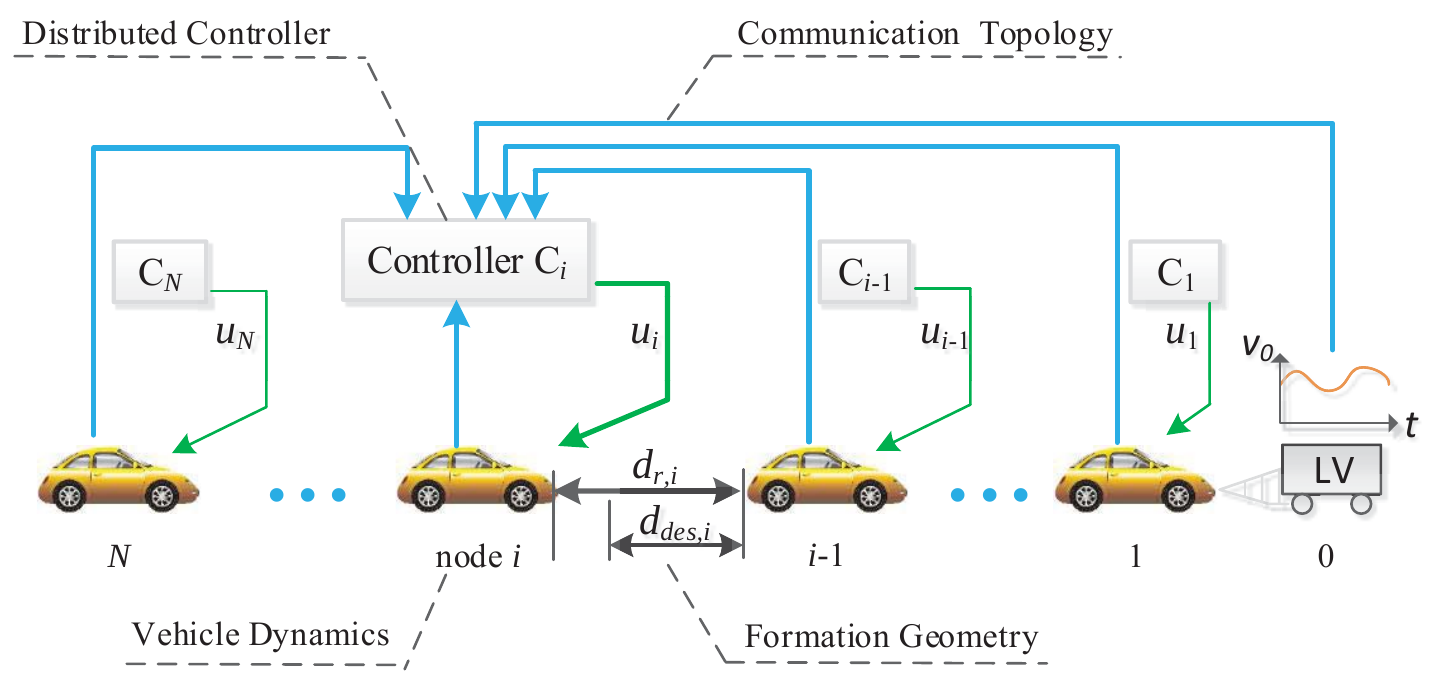}
    \caption{Four main components of a platoon \cite{zheng2015dynamic,zheng2015stabilityMargin}: 1) Vehicle dynamics; 2) Communication topology; 3) Distributed controller, 4) Formation geometry. $d_{r,i}$ denotes the actual relative distance; $d_{des,i}$ means the desired relative distance; $u_i$ is the control signal of $i$-th node; and $C_i$ represents the local controller in node $i$.}
    \label{fig:1}
\end{figure}

\section{Problem Statement} \label{section:pro}

This paper considers a homogeneous platoon of connected vehicles running on a flat road (see Fig. \ref{fig:1}), which has $N +1$ vehicles (or nodes), including a leading node indexed by $0$ and $N$ following nodes indexed from $1$ to $N$. The objective of platoon control is to ensure that the vehicles in a group to move at the same speed while maintaining a rigid formation geometry.

As demonstrated in Fig. \ref{fig:1}, a platoon can be viewed as a combination of four components: 1) vehicle dynamics; 2) communication topology; 3) distributed controller; and 4) formation geometry, which is called the four-component framework, originally proposed in \cite{zheng2015dynamic} and \cite{zheng2015stabilityMargin}. The vehicle dynamics corresponds to the behavior of each node; the communication topology describes information exchanging architecture among the nodes in a platoon; the distributed controller implements the feedback laws based on the available information for each node; and the formation geometry defines the desired inter-spacing between two adjacency nodes. Each component can exert different influence on the collective behavior of a platoon. Based on the features of each component, the existing literature on platoon control is categorized and summarized in \cite{zheng2015dynamic} and \cite{li2015overview}.

\subsection{Vehicle Longitudinal Dynamics}

The longitudinal dynamics of each vehicle is inherent nonlinear, including the engine, brake system, aerodynamics drag, rolling resistance, and gravitational force, \emph{etc}. However, a detailed nonlinear vehicle model may be not suitable for theoretical analysis. Many previous work either employed a hierarchical control framework consisting of a lower level controller and an upper level controller, or used a feedback linearization technique to get a linear vehicle model as a basis for theoretic analysis \cite{li2015overview}. Here, we use a first-order inertial function \eqref{eq:1} to approximate the acceleration response of vehicle longitudinal dynamics
\begin{equation}\label{eq:1}
  \tau\dot{a}_i(t) + a_i(t) = u_i(t)+ w_i(t),
\end{equation}
where $a_i(t)$ is the acceleration of node $i$; $\tau$ denotes the time delay in powertain systems; $u_i$ is the control input, representing the desired acceleration; and $w_i(t)$ denotes the external disturbance.

Choosing each vehicle's position $p_i(t)$, velocity $v_i(t)$ and acceleration $a_i(t)$ as the state, a state space representation of the vehicle dynamics is formulated as
\begin{equation}\label{eq:2}
  \dot{x}_i(t) = Ax_i(t)+ B_1u_i(t)+B_2w_i(t),
\end{equation}
where
\begin{equation*}
  x_i=\left[
        \begin{array}{c}
          p_i \\
          v_i \\
          a_i \\
        \end{array}
      \right]
  , A = \left[
        \begin{array}{ccc}
          0 & 1 & 0 \\
          0 & 0 & 1 \\
          0 & 0 & -\frac{1}{\tau} \\
        \end{array}
      \right],B_1=B_2=\left[
        \begin{array}{c}
          0 \\
          0 \\
          \frac{1}{\tau}  \\
        \end{array}
      \right] .
\end{equation*}

Note that the model \eqref{eq:2}, as well as its variants, is widely used as a basis of analysis in many vehicle control applications, \emph{e.g.}, \cite{naus2010string,oncu2014cooperative,zheng2015stabilityMargin, hao2013stability, wang2015longitudinal}.

\subsection{Model of Communication Topology}

Here, directed graphs are used to model the allowable communication connections between vehicles in a platoon. For more comprehensive descriptions on graph theory, please see \cite{godsil2013algebraic} and the references therein.  In this paper, it is assumed that the communication is perfect, and ignore the effects such as data quantization, time delay and switching effects.   % \cite{gao2016robust}

Specifically, we use a directed graph $\mathcal{G}_N = (\mathcal{V}_N,\mathcal{E}_N)$  to model the communication connections among the followers, with a set of $ N $ vertices $\mathcal{V}_N =\{1,2,\ldots,N\}$ and a set of edges $\mathcal{E}_N \subseteq \mathcal{V}_N \times \mathcal{V}_N$. We further associate an adjacency matrix $\mathcal{A} = [a_{ij}]_{N\times N}$ to the graph $ \mathcal{G} $: $a_{ij} = 1$ in the presence of a directed communication link from node $j$ to node $i$, \emph{i.e.}, $(j,i) \in \mathcal{E}_N$; otherwise $a_{ij} = 0$. Moreover, self-loops are not allowed here, \emph{i.e.}, $a_{ii} = 0$ for all $i\in \mathcal{V}_N$. The in-degree of node $i$ is defined as $d_i=\sum_{j=1}^Na_{ij}$. Denote $\mathcal{D}= \text{diag}\{d_1,d_2,\ldots,d_N\}$, and the Laplacian matrix $\mathcal{L}$ is defined as $\mathcal{L}=\mathcal{D}-\mathcal{A}$. We define a neighbor set of node $i$ in the followers as
\begin{equation}\label{eq:3}
  \mathbb{N}_i = \{j\in \mathcal{V}_N \mid a_{ij} = 1\}.
\end{equation}

To model the communications from the leader to the followers, we define an augmented directed graph $\mathcal{G}_{N+1} = (\mathcal{V}_{N+1},\mathcal{E}_{N+1})$ with a set of $ N+1 $ vertices $\mathcal{V}_{N+1} =\{0,1,\ldots,N\}$, which includes both the leader and the followers in a platoon. We use a pinning matrix $\mathcal{P}= \text{diag}\{p_1,p_2,\ldots,p_N\}$ to denote how each follower connects to the leader: $p_i = 1$ if $(0,i)\in \mathcal{E}_{N+1}$, otherwise $p_i = 0$. Note that $(0,i)\in \mathcal{E}_{N+1}$ means that node $i$ can obtain the leader's information via wireless communication, where we call node $i$ is pinned to the leader. The leader accessible set of node $i$ is defined as
\begin{equation}\label{eq:4}
  \mathbb{P}_i = \begin{cases}
                \{0\} & \text{  if } p_i = 1 \\
                \emptyset & \text{  if }p_i =0
                \end{cases}.
\end{equation}

This paper focus on undirected communication topology. The communication topology is called undirected if the information flow among followers (\emph{i.e.}, graph $\mathcal{G}_N$) is undirected, which means $i \in \mathbb{N}_j \Leftrightarrow j \in \mathbb{N}_i, \forall i,j \in \mathcal{V}_N$. Note that we do not restrict the number of followers that connects to the leader in an undirected topology. A spanning tree is a tree connecting all the nodes of a graph \cite{godsil2013algebraic}. Throughout this paper, we make the following assumption.

\begin{myAss}
    The augmented graph $\mathcal{G}_{N+1}$ contains at least one spanning tree rooting at the leader.
\end{myAss}

This assumption implies the leader is globally reachable in $\mathcal{G}_{N+1}$. In other words, every follower can obtain the leader information directly or indirectly, which is a prerequisite to guarantee the internal stability of a platoon.

\begin{myLem} \cite{olfati2004consensus}
    \label{lem:2}
    For any undirected communication topology, $\lambda_{\min}(\mathcal{L}) = 0$  with $\textbf{1}_N = [1,\ldots,1]^T \in \mathbb{R}^N$ as the corresponding eigenvector. Moreover, if Assumption 1 holds, $\lambda_{\min}(\mathcal{L}) = 0$ is a simple eigenvalue, and all the eigenvalues of $\mathcal{L} + \mathcal{P}$ are greater than zero, \emph{i.e.},  $\lambda_i(\mathcal{L} + \mathcal{P}) >0, i = 1, \ldots, N$.
\end{myLem}

Note that the aforementioned techniques are widely used in the consensus of multi-agent systems \cite{olfati2004consensus,li2011h,cao2013overview}, which have recently been applied to study the influence of different topologies on platoon performance in \cite{zheng2016stability,zheng2015dynamic} and \cite{zheng2015stabilityMargin}.

\subsection{Design of Linear Distributed Controller}

The objective of platoon control is to maintain the same speed with the leader and to keep a desired inter-vehicle space:%, \emph{i.e.},
\begin{equation}\label{eq:5}
  \begin{cases}
                \displaystyle{\lim_{t \to +\infty}}\|v_i-v_0\| = 0 \\%\lim\limits_{t\rightarrow +\infty }{} \\
                \displaystyle{\lim_{t \to + \infty}} \|p_i-p_{i-1}-d_{i,i-1}\|=0
  \end{cases}, i = 1,\ldots,N.
\end{equation}
where $d_{i,i-1}$ denotes the desired distance between node $i$ and node $i-1$, and $v_0$ is the leader's velocity. We assume that the leader runs a constant speed trajectory, \emph{i.e.}, $a_0 = 0, p_0 = v_0t$. In this paper, we use constant spacing policy, \emph{i.e.}, $d_{i,i-1}=d_0$, which is widely employed in the literature \cite{hao2011stability,barooah2009mistuning,zheng2015stabilityMargin,li2015overview}.

The synthesis of local controller for node $i$ can only use the information of nodes in set $\mathbb{I}_i=\mathbb{N}_i \cup \mathbb{P}_i$. Here, we use a linear identical state feedback form, as employed in \cite{hao2011stability,barooah2009mistuning,zheng2015stabilityMargin},
\begin{equation}\label{eq:6}
  u_i = -\sum_{j\in\mathbb{I}_i}c[k_p(p_i-p_{j}-d_{i,j})+k_v(v_i-v_j)+k_a(a_i-a_j)],
\end{equation}
where $k_p,k_v,k_a$ are the local feedback gains, and $c$ denotes the coupling strength which is convenient for the discussions in controller synthesis. In this paper, it is assumed the communication is perfect, and ignore the effects such as data quantization, time delay and switching effects.

To write a compact form of \eqref{eq:6}, we define the tracking errors for node $i$
\begin{equation}\label{eq:7}
  \left \{ \begin{array}{l}
                \hat{p}_i=p_i - p_0 - d_{i,0} \\%\lim\limits_{t\rightarrow +\infty }{} \\
                \hat{v}_i=v_i-v_0 \\
                \hat{a}_i=a_i -a_0
                \end{array}
              \right..
\end{equation}
Then, \eqref{eq:6} can be rewritten into
\begin{equation}\label{eq:8}
  u_i = -ck^T\sum_{j\in\mathbb{I}_i}(\hat{x}_i-\hat{x}_j), i = 1,\ldots,N,
\end{equation}
where $k=[k_p,k_v,k_a]^T$ is the vector of feedback gains, and $\hat{x}_i = [\hat{p}_i,\hat{v}_i,\hat{a}_i]^T$.

\begin{figure}[!t]
    \centering
    \includegraphics[width=0.95\columnwidth]{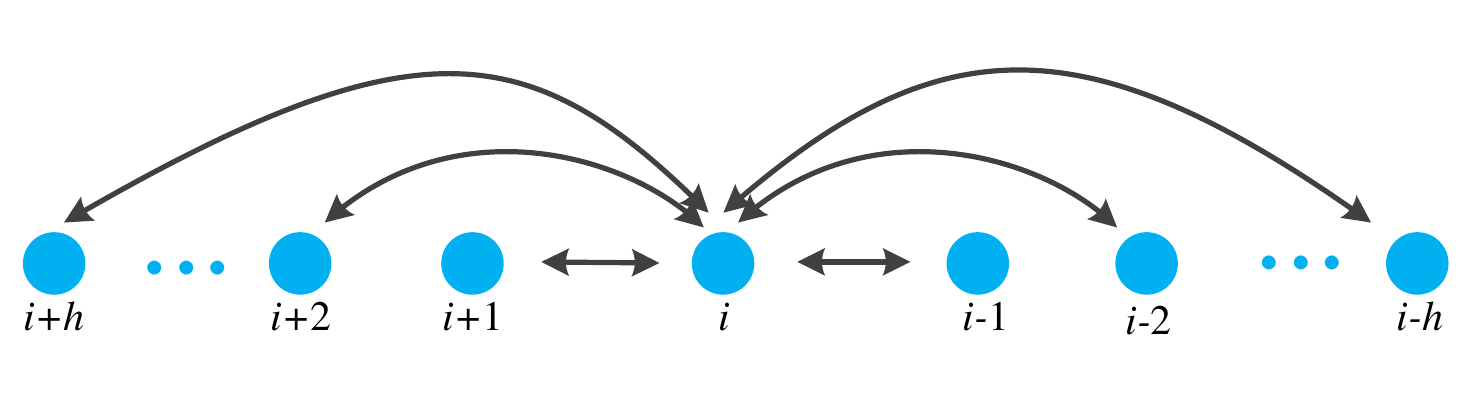}
    \caption{Illustration of $h$-neighbor undirected topology, where $h$ denotes the reliable communication range.}
    \label{fig:2}
\end{figure}

\subsection{Formulation of Closed-loop Platoon Dynamics}

The collective state and input of all following vehicles in a platoon are defined as $X=[\hat{x}_1^T,\hat{x}_2^T,\ldots, \hat{x}_N^T]^T$, and $U=[u_1,u_2,\ldots,u_N]^T$, respectively. Based on \eqref{eq:8}, we have
\begin{equation}\label{eq:9}
  U = -c(\mathcal{L}+\mathcal{P})\otimes k^T \cdot X,
\end{equation}
where $\otimes $ denotes the Kronecker product. Then, we have the closed-loop dynamics of the homogeneous platoon as
\begin{equation}\label{eq:10}
 \dot{X} = [I_N \otimes A - c(\mathcal{L}+\mathcal{P}) \otimes B_1k^T] \cdot X+B \cdot W,
\end{equation}
where $B= I_N \otimes B_2$, $I_N$ is the identity matrix of dimension $N\times N$, and $W = [w_1,w_2,\ldots,w_N]^T $ denotes the vector of external disturbances. Here, we denote
\begin{equation}\label{eq:Ac}
   A_c = I_N \otimes A - c(\mathcal{L}+\mathcal{P}) \otimes B_1k^T.
\end{equation}

We define the tracking error of positions as the output of a platoon, \emph{i.e.},
\begin{equation}\label{eq:11}
 Y = [\hat{p}_1,\hat{p}_2,\ldots,\hat{p}_N]^T = C \cdot X.
\end{equation}
where $C = I_N \otimes C_1$ and $ C_1 = [1,0,0]$. Assuming zero initial tracking errors, we obtain the transfer function from $W$ to $Y$ as
\begin{equation}\label{eq:12}
    G(s) = C(sI_{3N}-A_c)^{-1}B,
\end{equation}
where $s$ denotes the complex number frequency. Considering the expressions of $A_c, B, C$, we can further obtain that
\begin{multline}\label{eq:13}
    G(s) = [I_N \cdot (\tau s^3+s^2)+ \\ c(\mathcal{L}+\mathcal{P})\cdot (k_p+k_vs+k_as^2)]^{-1}.
\end{multline}

Note that \eqref{eq:10} and \eqref{eq:13} are unified models in time domain and frequency domain, respectively, which describe the platoon dynamics with various communication topologies. From \eqref{eq:13}, it is easy to find that the robustness of a platoon depends on not only the distributed feedback gains, but also the communication topology that interconnects the vehicles. Moreover, the communication topology can cast fundamental limitations on certain platoon properties, such as stability margin \cite{hao2011stability,zheng2015stabilityMargin}, string stability \cite{seiler2004disturbance}, and coherence behavior \cite{li2015overview,bamieh2012coherence}. In this paper, we focus on analyzing the scaling trend of the robustness and synthesizing the distributed controller with guaranteed performance for a large-scale platoon.

\section{Robustness Analysis of platoons with undirected topologies}\label{section:Robustness}

\begin{figure}[!t]
    \centering
    \includegraphics[width=0.89\columnwidth]{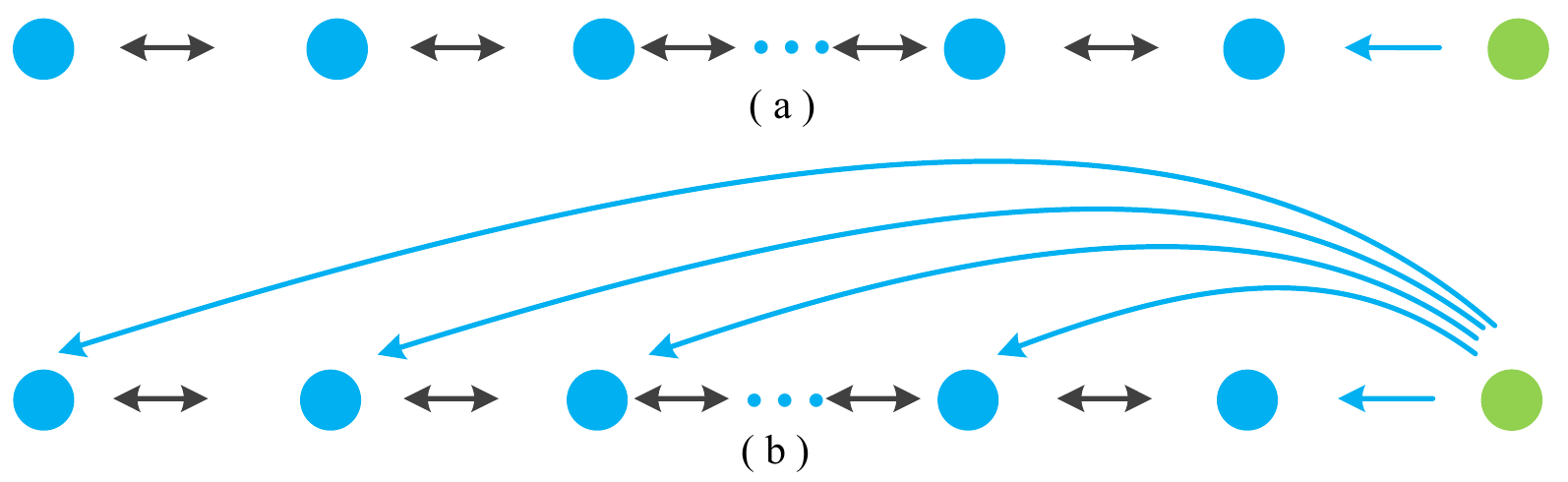}
    \caption{Typical examples of undirected topology. (a) Bidirectional (BD) topology; (b) Bidirectional-leader (BDL) topology}
    \label{fig:3}
\end{figure}

In this section, we present the results on robustness analysis of large-scale platoons with undirected communication topologies. In fact, an undirected topology includes a large class of topologies, where one commonly used variant is the $h$-neighbor undirected topology (see Fig. \ref{fig:2})~\cite{zheng2015stabilityMargin}.

\begin{myDef}
    ($h$-neighbor undirected topology) The communication topology is called to be an $h$-neighbor undirected topology, if each follower can reach its nearest $h$ neighbors in graph $\mathcal{G}_N$, \emph{i.e.}, $\mathbb{N}_i = \{i-h,\ldots,i+h\} \cap \mathcal{V}_N \setminus \{i\}$.
\end{myDef}

The parameter $h$ can be viewed as the reliable range of wireless communication. Typical examples of undirected topology are the bidirectional (BD) and bidirectional-leader (BDL) topology (where $h=1$; see Fig.~\ref{fig:3}(a) and Fig.~\ref{fig:3}(b) for illustration), which have been discussed in terms of string stability \cite{seiler2004disturbance} and stability margin \cite{zheng2015stabilityMargin}.

\begin{figure*}[hb]
    \hrulefill

    \newcounter{TempEqCnt}
    \setcounter{TempEqCnt}{\value{equation}}
    \setcounter{equation}{20}

    \begin{equation} \label{eq:31}
        \begin{split}
            \|G_i(s)\|_{\mathcal{H}_{\infty}} & = \sup_{\omega}\sqrt{\frac{1}{(\lambda_ik_v\omega - \tau \omega^3)^2+[\lambda_ik_p-(1+\lambda_ik_a)\omega^2]^2}} \\
            & = \sup_{\omega}\sqrt{\frac{1}{\tau^2\omega^6+[(1+\lambda_ik_a)^2-2\lambda_ik_v\tau]\omega^4+[(\lambda_ik_v)^2-2\lambda_ik_p(1+\lambda_ik_a)]\omega^2+(\lambda_ik_p)^2}} \\
            & \geq \frac{1}{\lambda_ik_p}
        \end{split} \quad .
    \end{equation}

    \setcounter{equation}{\value{TempEqCnt}}
\end{figure*}

\subsection{Robustness Measure: $\gamma$-gain}

Here, a performance measure is introduced for platoons with external disturbances from the perspective of energy amplification.
Specifically, we consider an appropriate amplification factor in the following scenario: disturbances acting on all vehicles $W \in \mathbb{R}^N$ to the position tracking errors of all vehicles $Y \in \mathbb{R}^N $\cite{hao2013stability}. In this paper, we only consider the disturbances with limited energy, which means the $\mathcal{L}_2$ norm of a disturbance $w_i$ is bounded, \emph{i.e.}, $ \|w_i(t)\|_{\mathcal{L}_2}=\int_0^{+\infty}|w_i(t)|^2 dt < \infty$. Considering the aforementioned scenario, we define the following $\gamma$-gain to quantify the robustness of a platoon.
\begin{myDef}
    ($\gamma$-gain) Consider a homogeneous platoon with the dynamics shown in \eqref{eq:10}. The $\gamma$-gain is defined as the following amplification:
        \begin{equation}\label{eq:15}
            \gamma = \sup \frac{\|Y(t)\|_{\mathcal{L}_2}}{\|W(t)\|_{\mathcal{L}_2}},
        \end{equation}
    where $Y(t)$ is the vector of tracking error of positions, and $W(t)$ is the vector of external disturbances. %with limited energy.
\end{myDef}

One physical interpretation is that $\gamma$-gain reflects the sensitivity or attenuation effect of the energy of external disturbances for a platoon. Note that the definition of $\gamma$-gain is not identical to the notion of standard string stability which reflects the attenuation effects of spacing error along the vehicle string \cite{seiler2004disturbance}. According to \cite{zhou1996robust}, $\gamma$-gain can be computed by using the $\mathcal{H}_{\infty}$ norms of corresponding transfer functions:
\begin{equation}\label{eq:17}
    \gamma = \|G(s)\|_{\mathcal{H}_{\infty}}= \sup_{\omega}[\sigma_{\text{max}} \big(G(j\omega)\big)],
\end{equation}
where $\sigma_{\text{max}}(\cdot)$ denotes the maximum singular value. % and $|\cdot|$ represents the complex modulus.

\begin{myRem}
    In principle, the feasibility to analyze the robustness of platoons depends on whether the $\mathcal{H}_{\infty}$ norms of certain transfer functions are explicitly computable. In general, it is rather difficult to analytically obtain the $\mathcal{H}_{\infty}$ norms for a general directed topology. We notice that one challenge comes from the inverse operation in transfer function $G(s)$ \eqref{eq:13}, especially when simultaneously considering the factor of communication topology $\mathcal{L}+\mathcal{P}$. Many previous works only focused on limited types of communication topologies, especially on predecessor-following and predecessor-following leader topology \cite{dunbar2012distributed,oncu2014cooperative,ploeg2014controller,seiler2004disturbance}, where the discussion on string stability is often case-by-case based. In this section, we explicitly analyze the scaling trend of the robustness for platoons with undirected topologies using the spectral decomposition of $\mathcal{L}+\mathcal{P}$.
\end{myRem}

\subsection{Scaling Trend of $\gamma$-gain for Large-scale Platoons with Undirected Topologies}

It is easy to know the matrix $\mathcal{L} + \mathcal{P}$ associated with any undirected topology is symmetric. For example, the matrix corresponding to bidirectional topology is listed as:
\begin{equation*}
    \mathcal{L}_{\text{BD}}+\mathcal{P}_{\text{BD}} = \left[
                                              \begin{array}{cccc}
                                                2 & -1 &  &  \\
                                                -1 & 2 & \ddots &  \\
                                                 &  \ddots & \ddots & -1 \\
                                                 &  & -1 & 1 \\
                                              \end{array}
                                            \right].
\end{equation*}

Before presenting the results on undirected topology, we need the following lemmas.

\begin{myLem} \label{lem:1}
    \cite{demmel1997applied} Given a symmetric matrix $Q \in \mathbb{S}^n$ and a real vector $x\in \mathbb{R}^{n}$, we have
    \begin{equation*}
        \begin{cases}
                \lambda_{\max}(Q) = \max \limits_{x \neq 0} \frac{x^TQx}{x^Tx} \\%\lim\limits_{t\rightarrow +\infty }{} \\
                \lambda_{\min}(Q) = \min \limits_{x \neq 0} \frac{x^TQx}{x^Tx}
       \end{cases}
    \end{equation*}
    where $\lambda_{\max}(Q), \lambda_{\min}(Q)$ are the maximum and minimum eigenvalue of matrix $Q$, respectively.
\end{myLem}

%Note that Lemma \ref{lem:1} is the well-known Rayleigh–Ritz theorem.

\begin{myLem} \cite{zheng2016stability} \label{lem:3}
    The minimum eigenvalue of $\mathcal{L}_{\text{BD}}+\mathcal{P}_{\text{BD}}$ satisfies
    \begin{equation*}
        \frac{1}{N^2} \leq  \lambda_{\min}(\mathcal{L}_{\text{BD}}+\mathcal{P}_{\text{BD}}) \leq \frac{\pi^2}{N^2}.
    \end{equation*}
\end{myLem}

Now we present the first result of this paper. Here, for robustness analysis, without loss of generality, we assume the coupling strength $c$ in the controller design is equal to one, \emph{i.e.}, $c = 1$ in \eqref{eq:6}.
\begin{myTheo} \label{theo:1}
    Consider a homogeneous platoon with undirected topology given by \eqref{eq:13}. Using any stabilizing feedback gains, the robustness measure $\gamma$-gain satisfies
    \begin{equation}\label{eq:28}
       \gamma \geq \frac {1}{\lambda_{\min}k_p},
    \end{equation}
    where $\lambda_{\min}$ is minimum eigenvalue of matrix $\mathcal{L} + \mathcal{P}$\footnote{In this paper, if not explicitly given the matrix, $\lambda_{\min}$ exclusively refers to the minimum eigenvalue of $\mathcal{L}+\mathcal{P}$.}.
\end{myTheo}

\begin{IEEEproof}
    Since $\mathcal{L} + \mathcal{P}$ is symmetric, there exists an orthogonal matrix $V \in \mathbb{R}^{N \times N}, VV^T = I_N$, such that
    \begin{equation} \label{eq:dec}
      \mathcal{L} + \mathcal{P} = V \Lambda V^T,
    \end{equation}
    where $\Lambda = \text{diag}\{\lambda_{1},\lambda_2, \ldots, \lambda_{N}\}$, and $\lambda_i$ is the $i$-th real eigenvalue of $\mathcal{L} + \mathcal{P}$.

    Then, based on \eqref{eq:13}, we know that the transfer function from disturbance $W$ to position tracking error $Y$ for platoons with undirected topology is
    \begin{equation*}
    \begin{split}
        G(s) & = [I_N (\tau s^3+s^2) + V \Lambda V^T (k_p+k_vs+k_as^2)]^{-1} \\
                         & = V[I_N (\tau s^3+s^2) + \Lambda (k_p+k_vs+k_as^2)]^{-1}V^T
    \end{split}.
    \end{equation*}
    Therefore, we have
    \begin{equation} \label{eq:29}
        G(s)  = V \left[
                                  \begin{array}{cccc}
                                    G_1 &  &  &\\
                                     & G_2 & &  \\
                                     &  & \ddots  &\\
                                     & & & G_N
                                  \end{array}
                                \right]V^T,
    \end{equation}
    where
    \begin{equation*}
        G_i(s) = \frac{1}{\tau s^3 + (1+\lambda_ik_a)s^2+\lambda_ik_vs+\lambda_ik_p}, i = 1,\ldots, N.
    \end{equation*}

    Thus, based on \eqref{eq:17}, we know that
    \begin{equation} \label{eq:30}
        \begin{split}
            \gamma & = \|G\|_{\mathcal{H}_{\infty}} = \sup_{\omega}\sqrt{\lambda[G^*(j\omega)G(j\omega)]} \\
            & = \sup_{\omega}\max_{i}\sqrt{\lambda[G^*_{i}(j\omega)G_{i}(j\omega)]} \\
            & = \max_i \|G_i(s)\|_{\mathcal{H}_{\infty}}
        \end{split}
    \end{equation}
    where $(\cdot)^{*}$ denotes the complex conjugate. Further, according to the expression of $G_i(s)$, we can explicitly calculate a lower bound of $\|G_i(s)\|_{\mathcal{H}_{\infty}}$, as shown in \eqref{eq:31}.

    Therefore, based on \eqref{eq:30}, \eqref{eq:31} and considering Lemma \ref{lem:2}, we have
    \setcounter{equation}{21}
    \begin{equation} \label{eq:32}
        \gamma \geq \frac{1}{\lambda_{\min}k_p}.
    \end{equation}
\end{IEEEproof}

Theorem \ref{theo:1} gives an explicit connection between the proposed robustness measure, \emph{i.e.} $\gamma$-gain, and underlying communication topology in a platoon. The key point of the proof lies in the spectral decomposition of a symmetric matrix, thus leading to the decomposition of the transfer function $G(s)$ (see \eqref{eq:29}). As such, we only need to analyze the $\mathcal{H}_{\infty}$ norm of certain transfer function corresponding to a single node, and this transfer function is modified by the eigenvalue of the communication topology. The spectral decomposition of symmetric matrices, as well as its variants, has been exploited in the consensus of multi-agent systems \cite{ren2008distributed,cao2013overview,olfati2004consensus}. More recently, chordal decomposition techniques have been utilized to decompose a large-scale system~\cite{zheng2017scalable}, which are promising to be applied in the analysis and synthesis of platoon systems.

\begin{myRem}
    According to \eqref{eq:32}, the index $\gamma$-gain is lower bounded by $\lambda_{\min}$, which agrees with the analysis of stability margin in~\cite{zheng2015stabilityMargin}.  Moreover, this result also gives hints to select the communication topology that associates with a larger minimum eigenvalue, in order to achieve better robustness performance for large-scale platoons. One practical choice is to increase the communication range that enlarges the minimum eigenvalue $\lambda_{\min}$~\cite{zheng2015stabilityMargin}.
\end{myRem}

Here, we further give two corollaries based on Theorem \ref{theo:1}.

\begin{myCor}\label{cor:1}
    Consider a homogeneous platoon with BD topology. Using any stabilizing feedback gains, the the robustness measure $\gamma$-gain satisfies
    \begin{equation} \label{eq:33}
        \gamma_{\text{BD}} \geq \frac{N^2}{k_p\pi^2}.
    \end{equation}
\end{myCor}

The proof is straightforward by applying Lemma \ref{lem:3} to \eqref{eq:28}.

\begin{myCor}\label{cor:2}
    Consider a homogeneous platoon with undirected topology given by \eqref{eq:13}. Using any stabilizing feedback gain, $\gamma$-gain increases at least as $\mathcal{O}(N)$, if the number of followers that are pinned to the leader is fixed.
\end{myCor}
\begin{IEEEproof}
    According to Lemma \ref{lem:1}, we have
    \begin{equation} \label{eq:34}
        \lambda_{\min} \leq \frac{x^T(\mathcal{L}+\mathcal{P})x}{x^Tx}, \forall x \neq 0.
    \end{equation}
    Then, by choosing $x = [1,1,\ldots,1]^T \in \mathbb{R}^{N \times 1}$, we know that $\mathcal{L}x = 0$ and
    \begin{equation*}
        x^T\mathcal{P}x = \sum_{i=1}^Np_i = \Omega(N).
    \end{equation*}
    where $\Omega(N)$ is the number of followers that are pinned to the leader. Therefore, combining Lemma \ref{lem:2} and \eqref{eq:34}, we have
    \begin{equation*}
        0 < \lambda_{\min} \leq \frac{\Omega(N)}{N}.
    \end{equation*}
    According to Theorem \ref{theo:1}, we know that
    \begin{equation*}
        \gamma \geq \frac{N}{\Omega(N)k_p}.
    \end{equation*}
    If $\Omega(N)$ is fixed and $k_p$ is given, then  the index $\gamma$-gain will increase at least as $\mathcal{O}(N)$.
\end{IEEEproof}

Both Corollaries \ref{cor:1} and \ref{cor:2} give a lower bound of robustness index $\gamma$-gain, implying that $\gamma$-gain increases as the growth of platoon size for certain undirected topology, regardless the design of feedback gains. Note that BD topology is a special type of undirected topology (where $\Omega(N) = 1$), and Corollary \ref{cor:1} is consistent with Corollary \ref{cor:2}, but gives a tighter bound.

Corollary \ref{cor:2} implies that the information from the leader plays a more important role than those among the followers. Intuitively, the communications among the followers help to regulate the local behavior. However, the objective of a large-scale platoon is to track the leader's trajectory. The information from the leader gives certain preview information of reference trajectory to followers, which is more important to guarantee the robustness performance. This hints us that the transmission of leader's information should have priority when the communication resources are limited in real implementations.

\section{Distributed $\mathcal{H}_{\infty}$ Controller synthesis of Vehicle Platoons} \label{section:controller}

In this section, we introduce a distributed synthesis method for platoons with guaranteed $\mathcal{H}_{\infty}$ performance: given a desired $\gamma_d$-gain, find the feedback gains ($k = [k_p,k_v,k_a]^T$) and coupling strength ($c$) such that $\gamma = \|G(s)\|_{\mathcal{H}_{\infty}} < \gamma_d$. We call this problem as the \emph{distributed $\mathcal{H}_{\infty}$ control problem} for vehicle platoons.

Here, we first present a decoupling technique that converts the distributed $\mathcal{H}_{\infty}$ control problem into a set of $\mathcal{H}_{\infty}$ control of independent systems that share the same dimension with a single vehicle. Then, a multi-step procedure is proposed for the $\mathcal{H}_{\infty}$ control problem of vehicle platoons, which inherits the favorable decoupling property. Finally, we further give several discussions on the selections of communication topology based on the analytical results.

\subsection{Decoupling of Platoon Dynamics}

The dimension of the closed-loop matrix that describes the collective behavior of a platoon \eqref{eq:10} is $3N \times 3N$, which is computationally demanding for large platoon size $N$ if we directly use dynamics \eqref{eq:10} for controller synthesis. Here, we employ the spectral decomposition of $\mathcal{L}+\mathcal{P}$ to decouple the platoon dynamics \eqref{eq:10}, as used in the robustness analysis in Section \ref{section:Robustness}, which offers significant computational benefits in the design of distributed $\mathcal{H}_{\infty}$ controller.

\begin{myLem} \label{lem:4}
    Given $\gamma_d > 0$, platoon \eqref{eq:10} is asymptotically stable and robustness measure $\gamma = \|G(s)\|_{\mathcal{H}_{\infty}} < \gamma_d$, if and only if the following $N$ subsystems are all asymptotically stable and the $\mathcal{H}_{\infty}$ norms of corresponding transfer functions are less than $\gamma_d$:
    \begin{equation}\label{eq:35}
      \begin{split}
        \dot{\bar{x}}_i &= (A-c\lambda_iB_1k^T)\bar{x}_i + B_2\bar{w}_i, \\
        \bar{y}_i &= C_1\bar{x}_i, i = 1, \ldots, N.
      \end{split}
    \end{equation}
\end{myLem}
\begin{IEEEproof}
    According to \eqref{eq:30}, we have
    \begin{equation*}
     \begin{aligned}
      \|G(s)\|_{\mathcal{H}_{\infty}} < \gamma_d
      \Leftrightarrow \max_i \|G_i(s)\|_{\mathcal{H}_{\infty}} < \gamma_d.
      \end{aligned}
    \end{equation*}
    It is easy to check that \eqref{eq:35} is a state-space representation of $G_i(s)$ \eqref{eq:29} (where $c = 1$).

    This observation completes the proof.
\end{IEEEproof}

Lemma \ref{lem:4} decouples the collective behavior of a platoon~\eqref{eq:10} into a set of $N$ individual subsystems \eqref{eq:35}.    Similar to the frequency domain, the influence of communication topology reflects on the fact that the decoupled systems are modified by the eigenvalues of $\mathcal{L}+\mathcal{P}$.  Also, the distributed $\mathcal{H}_{\infty}$ control problem of vehicle platoons \eqref{eq:10} is equivalently decomposed into a set of $\mathcal{H}_{\infty}$ control problems of subsystems sharing the same dimension with a single vehicle in \eqref{eq:2}. Note that the dimension of each decoupled system in \eqref{eq:35} is only $3 \times 3$, which leads to a significant reduction in terms of computational complexity for controller synthesis.

\subsection{Synthesis of Distributed $\mathcal{H}_{\infty}$ Controller}

We first present one useful lemma, which provides a necessary and sufficient condition to grantee the existence of a controller \eqref{eq:6} for the distributed $\mathcal{H}_{\infty}$ control problem.

\begin{myLem} \label{lem:7} \cite{li2011h}
     Consider a linear time-invariant system described by \eqref{eq:10}. For a desired $\gamma_d > 0$, there exists a controller \eqref{eq:6} such that $\gamma = \|G(s)\|_{\mathcal{H}_{\infty}} < \gamma_d$, if and only if there exist a matrix $Q \succ 0$ and a scalar $\alpha>0$ such that
    \begin{equation} \label{eq:LMI}
      \left[
        \begin{array}{ccc}
          AQ+QA^T-\alpha B_1B_1^T & B_2 & QC_1^T\\
          B_2^T & -\gamma^2_dI & 0\\
          C_1Q & 0 & -I\\
        \end{array}
      \right] \prec 0.
    \end{equation}
\end{myLem}

Next, we introduce the second theorem of this paper.

\begin{myTheo} \label{theo:2}
   Consider a homogeneous platoon with undirected topology given by \eqref{eq:10}. For any desired $\gamma_d > 0$, we have the following two statements.
   \begin{enumerate}
     \item There always exists a feasible solution: $Q \succ 0$ and $\alpha >0$ to \eqref{eq:LMI};
     \item For the closed-loop platoon system, we have $|G(s)\|_{\mathcal{H}_{\infty}} < \gamma_d$, if choosing the feedback gains as $k^T = \frac{1}{2} B_1^TQ^{-1}$, and selecting the coupling strength satisfying
         \begin{equation}\label{eq:47}
            c \geq \frac{\alpha}{\lambda_{\min}}.
        \end{equation}
        where $Q \succ 0$ and $\alpha >0$ is a feasible solution to \eqref{eq:LMI}.
   \end{enumerate}
\end{myTheo}

\begin{IEEEproof}
    1) Using the Schur complement, we know \eqref{eq:LMI} is equivalent to
     \begin{equation}\label{E:Schur}
            AQ+QA^T - \alpha B_1B_1^T + \frac{1}{\gamma_d^2}B_2B_2^T+QC_1^TC_1Q \prec 0.
    \end{equation}
    According the vehicle dynamics \eqref{eq:2}, we have $B_1=B_2$. Then \eqref{E:Schur} is reduced to
    \begin{equation*}
            AQ+QA^T - (\alpha- \frac{1}{\gamma_d^2})B_1B_1^T+QC_1^TC_1Q \prec 0.
    \end{equation*}
    Therefore, it is easy to know that for any small $\gamma_d >0$, there always exists a feasible solution: $Q \succ 0$ and $\alpha >0$ to \eqref{eq:LMI} (by just increasing $\alpha$).

    2) According to Lemma \ref{lem:4}, $\|G(s)\|_{\mathcal{H}_{\infty}} < \gamma_d$ if there exists some feedback gains $k$, and coupling strength $c$, such that
    \begin{equation}\label{eq:40}
      \begin{cases}
                A-c\lambda_iB_1k^T \text{ is Hurwitz} \\%\lim\limits_{t\rightarrow +\infty }{} \\
                \|C_1(A-c\lambda_iB_1k^T)^{-1}B_2\|_{\mathcal{H}_{\infty}} < \gamma_d
       \end{cases} i = 1,\ldots, N.
    \end{equation}

    For notional simplicity, we let $A_{c_i} = A-c\lambda_iB_1k^T$. Considering the expression of controller $k^T = \frac{1}{2} B_1^TQ^{-1}$, we get
    \begin{equation} \label{eq:45}
       A_{c_i}Q+QA_{c_i}^T = AQ+QA^T - c\lambda_iB_1B_1^T.
    \end{equation}
    If the coupling strength satisfies \eqref{eq:47}, we know
    $$c\lambda_i \geq \alpha, i= 1, \ldots, N.  $$
    Then, considering \eqref{E:Schur} and \eqref{eq:45}, we have
    \begin{equation} \label{eq:46}
       A_{c_i}Q+QA_{c_i}^T + \frac{1}{\gamma_d^2}B_2B_2^T  + QC_1^TC_1Q \prec 0,
    \end{equation}
    According to the Bounded Real Lemma \cite{zhou1996robust}, \eqref{eq:46} leads to \eqref{eq:40}. This completes the proof.

\end{IEEEproof}

To synthesize a distributed $\mathcal{H}_{\infty} $ controller, we only need to solve LMI \eqref{eq:LMI} to get the feedback gains $k^T$, and adjust the coupling strength $c$ to satisfy condition \eqref{eq:47}. Note that the dimension of LMI \eqref{eq:LMI} only corresponds to the dynamics of one vehicle, which is independent of the platoon size $N$. This fact reduces the computational complexity significantly, especially for large-scale platoons. This is a scalable multi-step procedure to design a distributed $\mathcal{H}_{\infty} $ controller for a platoon system. %The generalization of this approach to multi-agent systems can be found in \cite{li2011h}.

\begin{myRem}
    The synthesis of distributed $\mathcal{H}_{\infty} $ controller is actually decoupled from the design of communication topology. Specifically, the feasibility of LMI \eqref{eq:LMI} is independent of the communication topology, and the communication topology only exerts influence on the controller \eqref{eq:6} through the condition \eqref{eq:47}. This implies that we may separate the controller synthesis and the topology selection into two different stages for the design of a platoon system.
\end{myRem}

\begin{myRem} \label{re:10}
     The first statement of Theorem \ref{theo:2} also implies that we can design a distributed controller that satisfies any given $\mathcal{H}_{\infty}$ performance for a platoon system. The price is that the coupling strength $c$ becomes very large (see \eqref{eq:47}), thus resulting in a high-gain controller, which is not practical in real implementations considering the saturation of actuators.
\end{myRem}

\subsection{Selection of Communication Topology} \label{section:communication}

Here, we present some discussions on the selection of communication topology for vehicle platoons based on the aforementioned analytical results. In this framework, the influence of communication topology on the collective behavior of a platoon is represented by the eigenvalues of $\mathcal{L}+\mathcal{P}$. Especially, the least eigenvalue $\lambda_{\min}$ has significant influence on both the scaling trend of robustness measure $\gamma$-gain (see \eqref{eq:28}) and the synthesis of distributed $\mathcal{H}_{\infty} $ controller (see \eqref{eq:47}).

Based on \eqref{eq:28} and \eqref{eq:47}, it is favorable that the matrix $\mathcal{L}+\mathcal{P}$ associated with the communication topology has a larger minimum eigenvalue $\lambda_{\min}$, which not only improves the robustness for a given controller, but also reduces the value of feedback gains for a given $\mathcal{H}_{\infty}$ performance. In a practical application, the communication resources might be limited. Here, if we assume the number of communication links is limited, then a good undirected communication topology can be obtained from the following optimization problem.
\begin{equation}\label{eq:48}
  \begin{split}
        & \qquad \max \lambda_{\min} \\
        \text{s.t.  }  & \begin{cases} \|\mathcal{G}_{N+1}\|_0 \leq N_c \\ \text{Assumption 1 holds} \end{cases},
  \end{split}
\end{equation}
where $\lambda_{\min}$ is the minimal eigenvalue of $\mathcal{L}+\mathcal{P}$, $\|\mathcal{G}_{N+1}\|_0$ denotes the number of communication links, and $N_c$ represents the communication resources.

\begin{myRem}
    Problem \eqref{eq:48} is only an intuitive description. The exact mathematical formulation, as well as its corresponding solutions, is one of the future work. Also, this optimization problem should consider the issue of local safety, such as the ability to avoid collisions between consecutive vehicles, in the design of a platoon system.
\end{myRem}

\begin{figure}[!t]
    \centering
    \includegraphics[width=0.9\columnwidth]{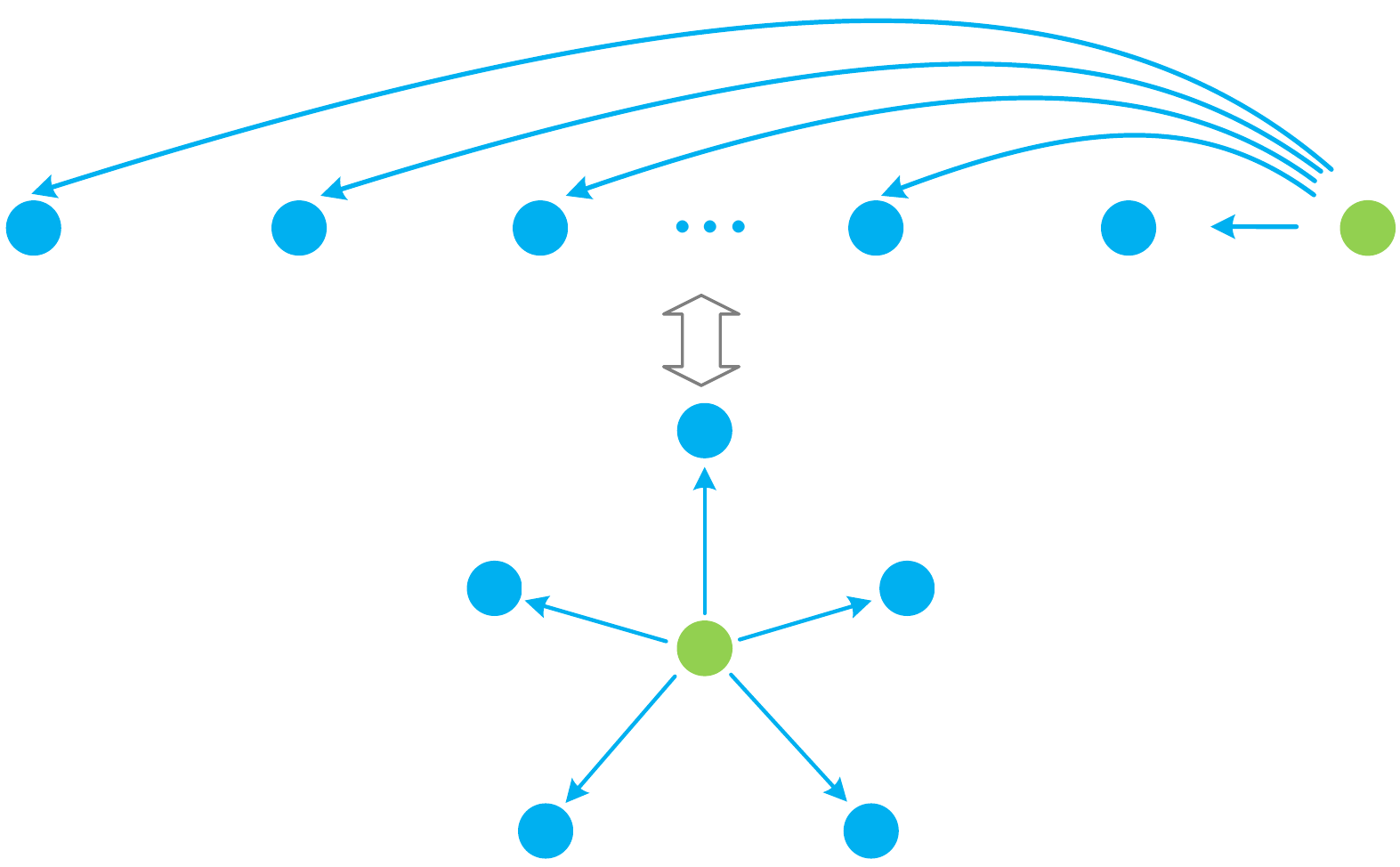}
    \caption{Communication topology with largest $\lambda_{\min}$ and minimum number of communication links: Star topology, where $\mathcal{P} = I$ and $\mathcal{L}=0$.}
    \label{fig:4}
\end{figure}

Here, we have the following results.

\begin{myLem} \label{lem:6} \cite{horn2012matrix}
    Given two symmetric matrices $A,B \in \mathbb{S}^{n}$, we have
    \begin{equation} \label{eq:lem6}
      \lambda_{\min}(A+B) \leq \lambda_{\min}(A)+\lambda_{\max}(B),
    \end{equation}
    where the equality holds if and only if there is a nonzero vector $x$ such that $Ax = \lambda_{\min}(A)x, Bx = \lambda_{\max}(B)x$, and $(A + B)x = \lambda_{\min}(A+B)x$.
\end{myLem}

\begin{myTheo} \label{theo:4}
    Given an undirected communication topology satisfying Assumption 1. The following statements hold:
    \begin{enumerate}
      \item $0 < \lambda_{\min} \leq 1$;
      \item $\lambda_{\min} = 1$ if and only if $\mathcal{P} = I$, \emph{i.e.}, every follower is pinned to the leader.
    \end{enumerate}
\end{myTheo}
\begin{IEEEproof}
    According to Lemma \ref{lem:2}, we know $ \lambda_{\min}>0$. Therefore, we only need to show $ \lambda_{\min}\leq1$.

    By Lemma \ref{lem:2}, $\lambda_{\min}(\mathcal{L}) = 0$. Further, based on the definition of $\mathcal{P}$, it is easy to know, $\lambda_{\max}(\mathcal{P}) = 1$. Therefore, according to Lemma \ref{lem:6}, we have
    \begin{equation}\label{eq:49}
        \lambda_{\min}(\mathcal{L}+\mathcal{P}) \leq \lambda_{\min}(\mathcal{L}) + \lambda_{\max}(\mathcal{P}) = 1.
    \end{equation}
    Thus, the first statement holds.

    Considering the equality condition in Lemma \ref{lem:6}, we have
    \begin{equation}\label{eq:50}
      \lambda_{\min} = 1 \Leftrightarrow \begin{cases} \mathcal{L}x=\lambda_{\min}(\mathcal{L})x = 0 \\ \mathcal{P}x=\lambda_{\max}(\mathcal{P})x = x \\ (\mathcal{L}+\mathcal{P})x = x \end{cases} ,x \neq 0.
    \end{equation}

    By lemma \ref{lem:2}, we have $\mathcal{L}x = 0~(x \neq 0) \Leftrightarrow x = \textbf{1}_N $. Considering the fact that $\mathcal{P}$ is diagonal, we know that
    \begin{equation}\label{eq:50}
        \mathcal{P}\textbf{1}_N = \textbf{1}_N \Leftrightarrow \mathcal{P}=I.
    \end{equation}

\end{IEEEproof}

Theorem \ref{theo:4} gives explicit upper and lower bounds for $\lambda_{\min}$, which are also tight. Besides, $\lambda_{\min}$ reaches the largest value when all the followers are pinned to the leader, regardless of the communication connections among followers. This fact also shows the importance of the leader's information, which is consistent with the robustness analysis in Section \ref{section:Robustness}.

\begin{figure}[!t]
    \centering
    \includegraphics[width=0.91\columnwidth]{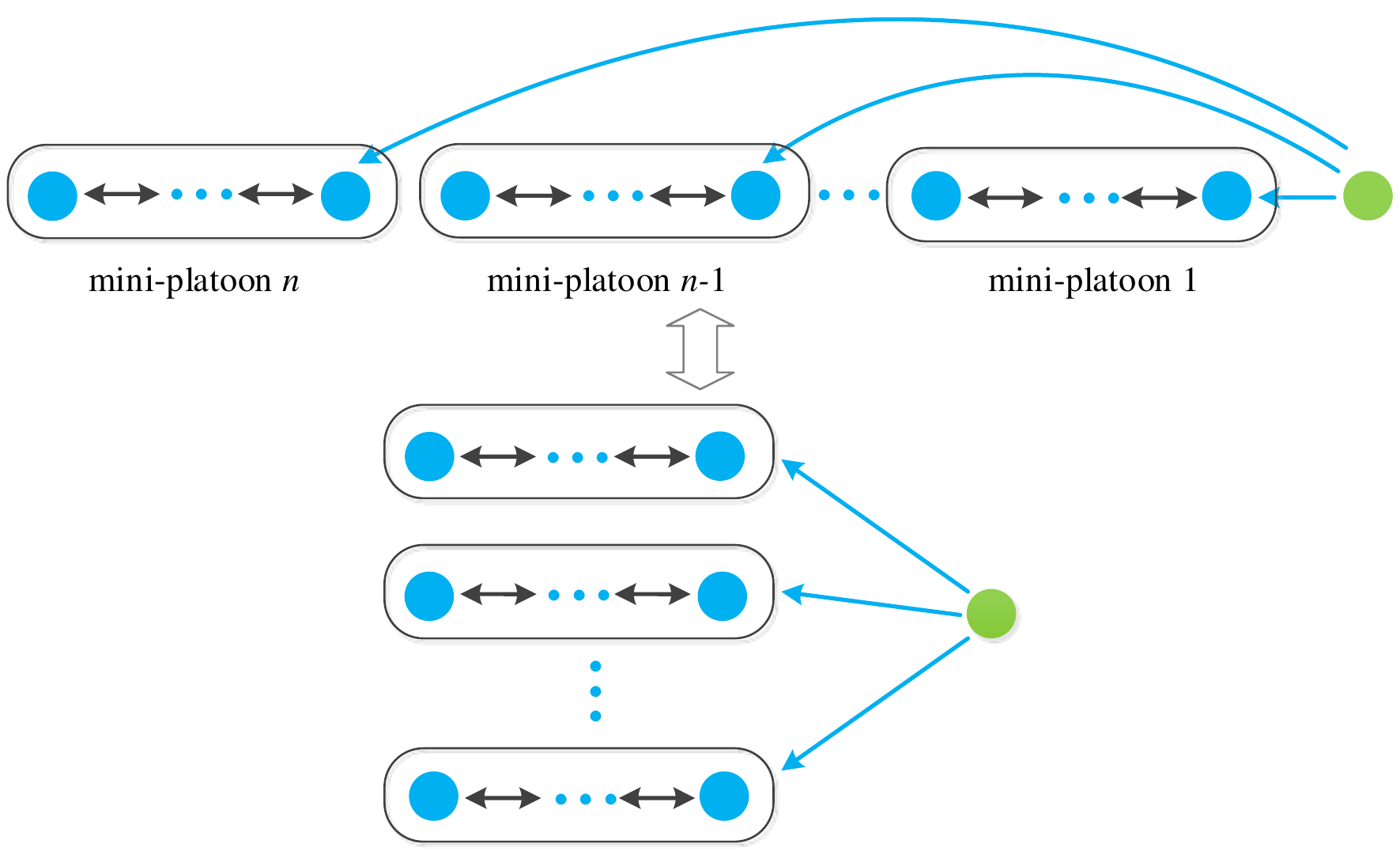}
    \caption{Coordination of multiple mini-platoons, where only the first node of each mini-platoon is pinned to the leader.}
    \label{fig:5}
\end{figure}

If we only limit the number of communication links in $\mathcal{G}_{N+1}$, the \emph{`best'} communication topology should be a star topology (see Fig. \ref{fig:4}), where all the followers are pinned to the leader and there exist no communications between followers, \emph{i.e.}, $\mathcal{P}=I$ and $\mathcal{L}=0$. In this case, $\lambda_{\min}$ reaches the largest possible value, and the number of communication links is minimum under Assumption 1. Here, note that we only consider the global behavior (\emph{i.e.}, reaching consensus with the leader), and ignore the local behavior (\emph{e.g.}, collisions between nearest vehicles). In a practical platoon, however, the vehicles' control decisions should take into account the behavior of its nearest neighbors. In this sense, BDL topology (see Fig. \ref{fig:3}(b)) might be a good choice.

\begin{myRem}
    In the discussions above, we assume that the communication is perfect. In reality, there always exist certain time-delays or data-losses in the communication channels, especially considering long-range and high-volume transmission of leader's data. Thus, it might not be very practical to let all the followers pinned to the leader for a large platoon. A balanced choice is to divide a platoon into several mini-platoons, where the first node of each mini-platoon is pinned to the leader, thus reducing the communication demanding from the leader (see Fig. \ref{fig:5} for illustrations). This kind of strategy can be viewed as the coordination of multiple mini-platoons~\cite{swaroop1999constant}, which is a higher level of platoon control.
\end{myRem}

\section{Numerical Simulations}\label{section:simulation}

In this section, numerical experiments with a platoon of passenger cars are conducted to verify our findings. We demonstrate both the scaling trend of $\gamma$-gain under different undirected topologies, and the computation of a distributed $\mathcal{H}_{\infty}$ controller. In addition, simulations with a realistic model of nonlinear vehicle dynamics (as used in~\cite{wang2015longitudinal,Zheng2016distributed}) are carried out to show the effectiveness of our approach in real environments.

\subsection{Scaling Trend of $\gamma$-gain}

In this subsection, we set the inertial delay as $\tau = 0.5$ s, and distributed feedback gains are chosen as $k_p = 1, k_v = 2, k_a = 0.5$, which stabilizes the platoons considered in this paper (see, \emph{e.g.}, the discussions on stability region of linear platoons in \cite{zheng2015stabilityMargin,zheng2016stability}).

\begin{figure}[!t]
    \centering
    \includegraphics[width=0.60\columnwidth]{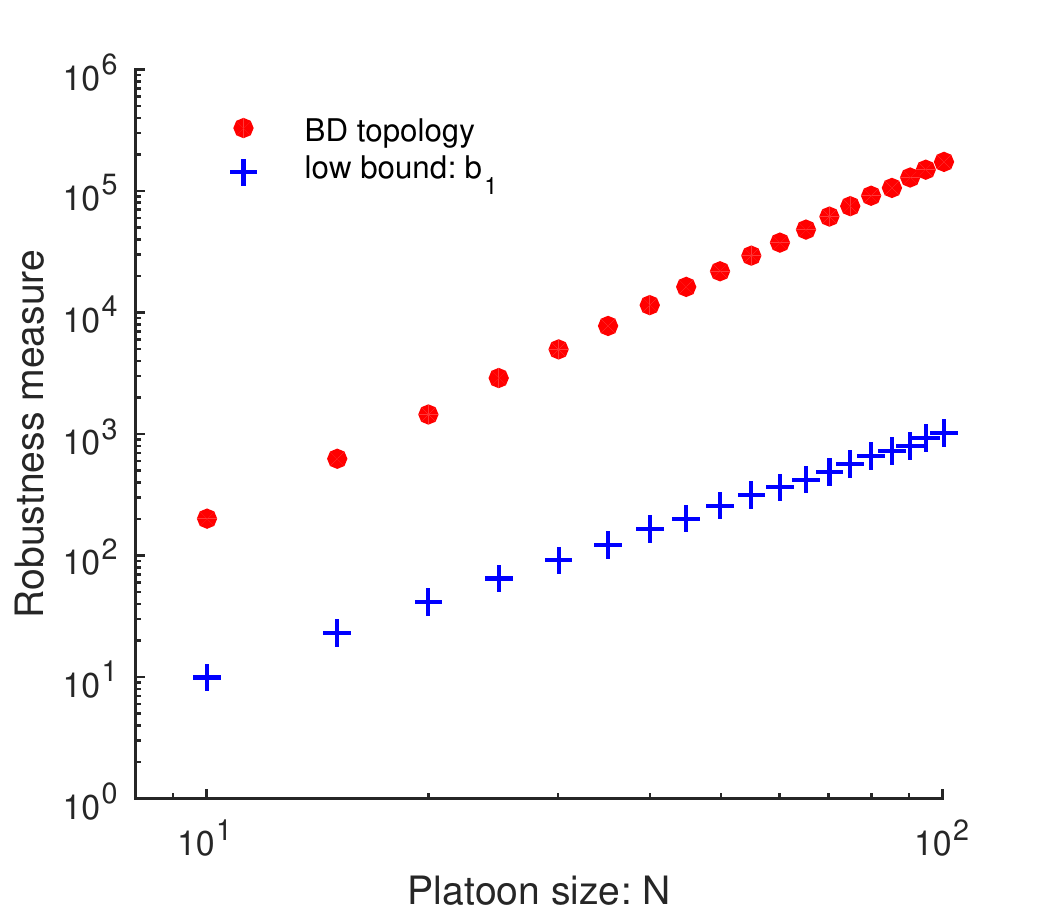}
    \caption{Scaling trend of $\gamma$-gain for platoons under BD topology (see Fig. \ref{fig:3}): lower bound $ b_1 = \frac{N^2}{k_p\pi^2}$ (see Corollary \ref{cor:1}). }
    \label{fig:6}
\end{figure}
\begin{figure}[!t]
    \centering
    \includegraphics[width=0.60\columnwidth]{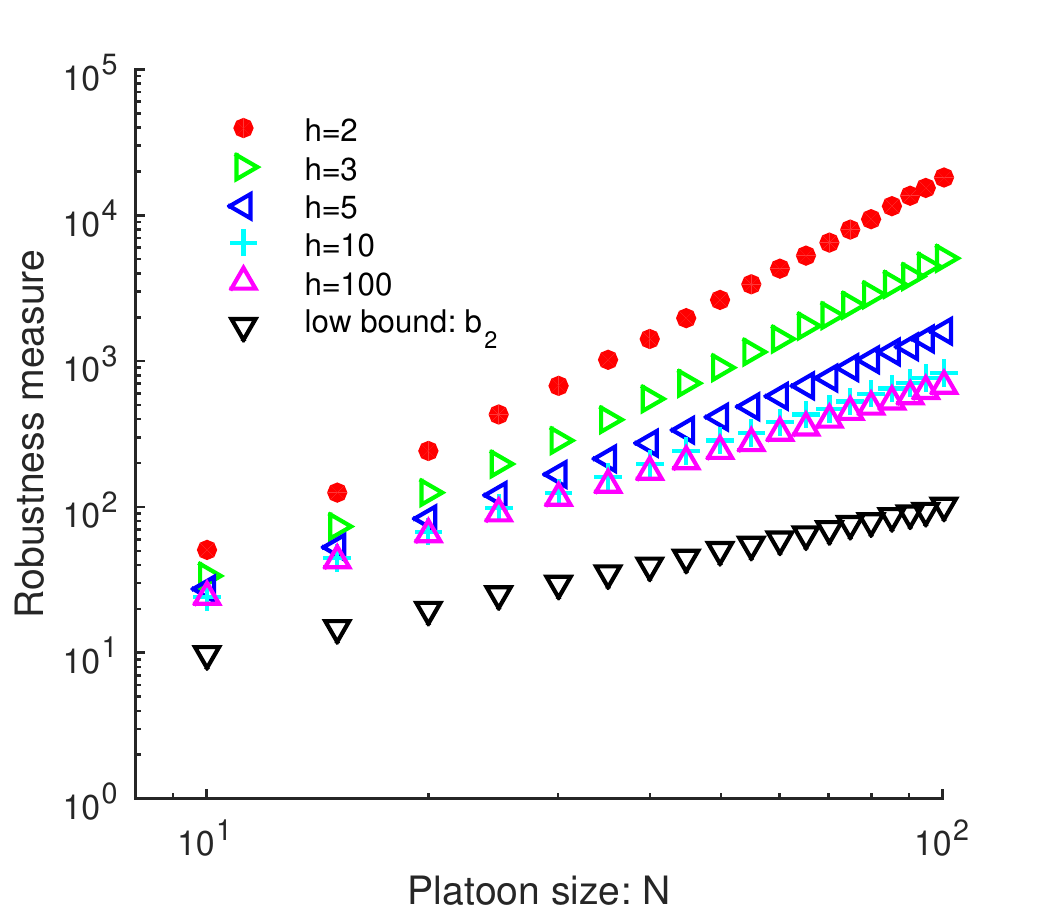}
    \caption{Scaling trend of $\gamma$-gain for platoons under undirected topology with different reliable communication range $h$ (see Fig. \ref{fig:2}): lower bound $ b_2 = \frac{N}{k_p\Omega(N)}$ (see Corollary \ref{cor:2}). Here, only the first follower is pinned to the leader (\emph{i.e.}, $\Omega(N) = 1$).}
    \label{fig:7}
\end{figure}
To illustrate the scaling trend of the robustness performance for platoons, we vary the platoon size $N$ from $10$ to $100$, and numerically calculate the $\mathcal{H}_{\infty}$ norm (\emph{i.e.}, $\gamma$-gain of this paper) and the lower bound of corresponding transfer functions. Fig. \ref{fig:6} demonstrates the scaling trend of $\gamma$-gain for BD topology, which indeed has a polynomial growth with the increase of platoon size (at least $\mathcal{O}(N^2)$). We can also find that the lower bound in Corollary \ref{cor:1} is mathematically correct but not very tight. How to get a tight lower bound and upper bound is one of our future work.

In a BD topology, the reliable communication range $h$ equal to one, which means each follower can only get access to the information of its direct front and back vehicle. Now, we examine the influence of different $h$ on the scaling trend of $\gamma$-gain. As shown in Fig. \ref{fig:7}, the increase of $h$ would slightly improve the robustness performance, but cannot alternate the scaling trend with the increase of platoon size (at least $\mathcal{O}(N)$), which obviously confirms with the statement in Corollary \ref{cor:2}. Note that when $h = 100$, the followers are fully connected in Fig. \ref{fig:7}, which means each follower can have full information of all other followers. In this case, the $\gamma$-gain still increase as the growth of platoon size, which agrees with the statements of the robustness analysis in Section III.

\subsection{Calculation of Distributed $\mathcal{H}_{\infty}$ Controller} \label{section:SimulationII}

\begin{figure}
    \centering
    \subfigure[ ]
    { \label{fig:a}
        \includegraphics[width=0.44\columnwidth]{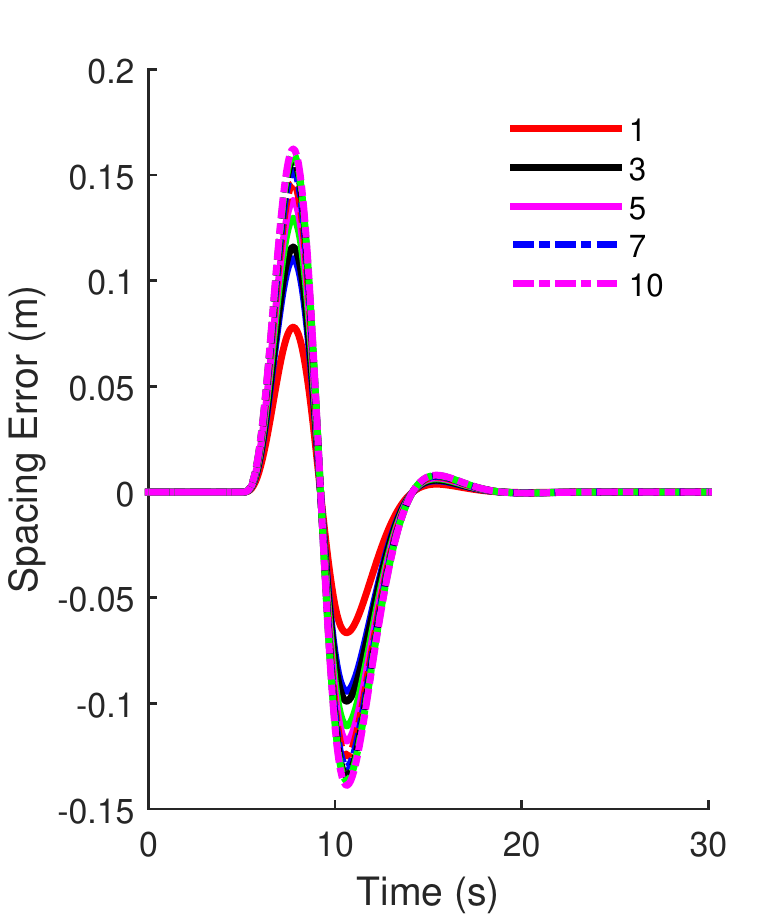}
    }
    \subfigure[ ]
    { \label{fig:b}
      \includegraphics[width=0.44\columnwidth]{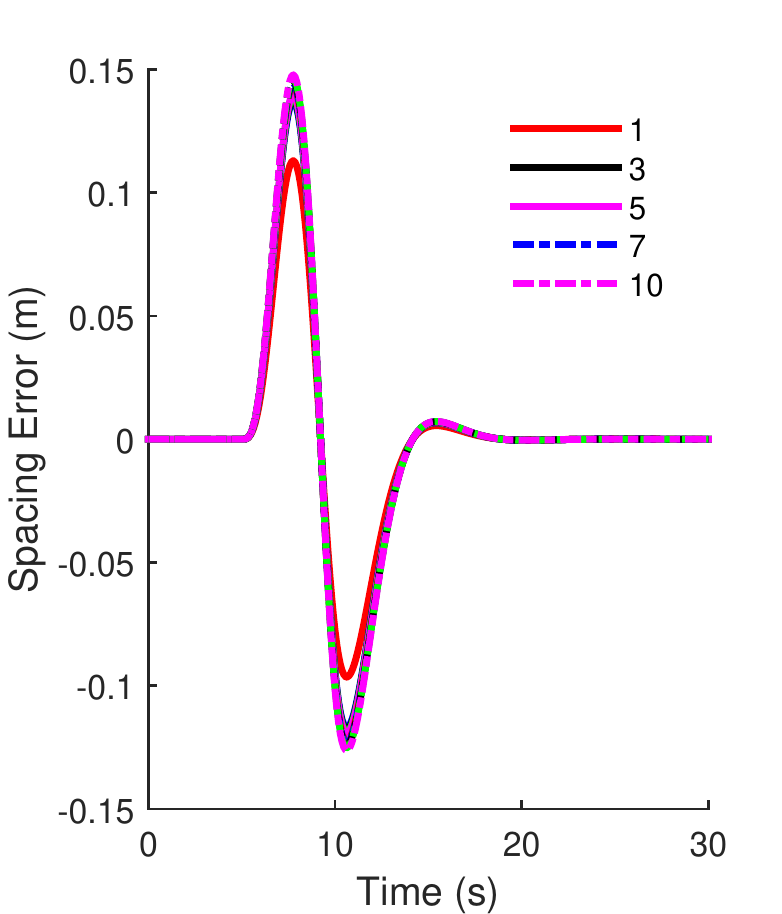}
    }
    \vspace{-0.4mm}
     \subfigure[ ]
    { \label{fig:a}
        \includegraphics[width=0.44\columnwidth]{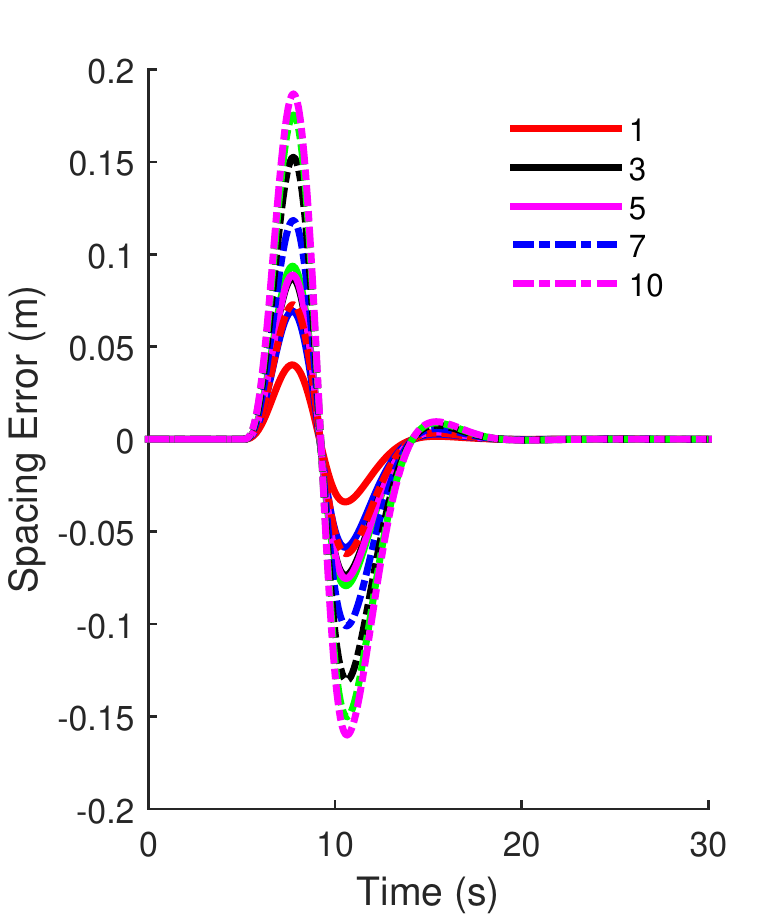}
    }
    \subfigure[ ]
    { \label{fig:b}
      \includegraphics[width=0.44\columnwidth]{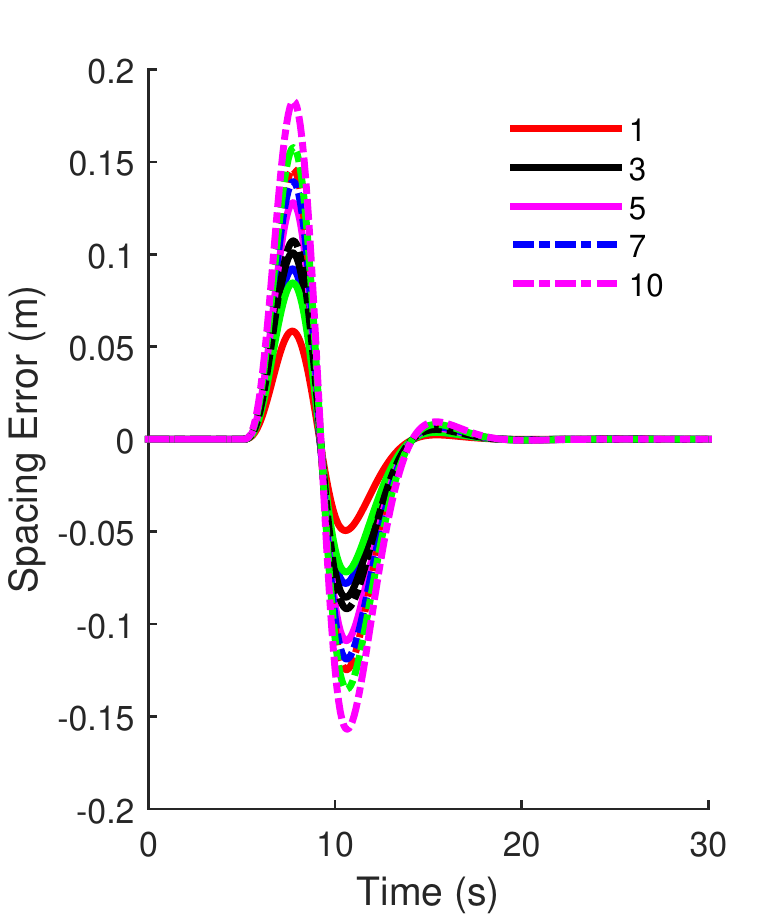}
    }
    \caption{Spacing errors for platoons with distributed $\mathcal{H}_{\infty}$ controllers (total platoon size $N=10$): (a) $h$-neighbor undirected topology with $h = 2$; (b) $h$-neighbor undirected topology with $h = 4$; (c) two multiple mini-platoons (both sizes are 5); (d) three multiple mini-platoons (the sizes are 3, 4 and 3).}
    \label{fig:9}
\end{figure}

The scalable multi-step procedure to synthesize a distributed $\mathcal{H}_{\infty}$ is summarized in Theorem \ref{theo:2}. As an example, we assume the inertial delay as $\tau = 0.5$ s. Solving the LMI \eqref{eq:LMI} with a desired performance $\gamma_d = 1$ using toolboxes YALMIP \cite{lofberg2004yalmip} and solver SeDuMi \cite{sturm1999using}, we get a feasible solution:
\begin{equation*}
  Q = \left[
        \begin{array}{ccc}
          0.669 & -0.419 & 0.006 \\
          -0.419 & 0.606 & -0.474 \\
          0.006 & -0.474 & 1.044 \\
        \end{array}
      \right], \alpha = 1.968.
\end{equation*}
Thus, according to Thoerem \ref{theo:2}, the feedback gain matrix is chosen as
\begin{equation} \label{eq:NumeController}
  k^T = [2.122 \quad 3.425 \quad 2.501].
\end{equation}

Then, for different communication topologies, we need to choose different coupling strength $c$ according to \eqref{eq:47}. For example, we consider a platoon of 11 vehicles (one leader and ten followers) with four types of topologies (see Fig. \ref{fig:2} and Fig. \ref{fig:5}), specified as follows:
\begin{itemize}
  \item a) $h$-neighbor undirected topology ($h = 2$);
  \item b) $h$-neighbor undirected topology ($h = 4$);
  \item c) two multiple mini-platoons (both size are 5);
  \item d) three multiple mini-platoons (the sizes are 3, 4 and 3).
\end{itemize}
Finally, we can choose the coupling strength according to \eqref{eq:47} for these four topologies, as listed in Table \ref{tab:1}. To achieve the desired performance $\gamma_d = 1$, Theorem \ref{theo:2} results in high-gain controllers, due to the low value of $\lambda_{\min}$ (see Remark \ref{re:10}).

We implement these controllers for the aforementioned platoons considering the following scenario: the initial state errors of the platoon are zeros, the leader runs a constant speed trajectory ($v_0 = 20~m/s$), and there exist the following external disturbances for each followers:% (see \eqref{eq:2}):
\begin{equation*}
  w_i(t) = \begin{cases}
                0 &0 < t <5 s \\
                \sin(\frac{2\pi}{5}(t-5)) & 5 s \leq t<10 s \\
                0 &  t \geq 10 s
                \end{cases}.
\end{equation*}
\begin{table}
    \centering
    \renewcommand\arraystretch{1.}
    \caption{Calculating coupling strength}
    \label{tab:1}
    \begin{tabular}{ c m{1.8cm}<{\centering} m{1.2cm}<{\centering} }
        \hline \toprule[1pt]
        Topologies & $\lambda_{\min}(\mathcal{L}+\mathcal{P})$ & $c$  \\    \hline
        a) & 0.0557 & 35.33 \\
        b) & 0.0806 & 24.42 \\
        c) & 0.0810 & 24.30 \\
        d) & 0.1790 & 10.99 \\
        \bottomrule[1pt]
        \end{tabular}
\end{table}

Fig. \ref{fig:9} shows the profiles of spacing errors in time-domain for these four types of homogeneous platoons, which are obviously stable under the $\mathcal{H}_{\infty}$ controller calculated by Theorem \ref{theo:2}. Further, we can calculate the error amplification $\frac{\|Y(t)\|_{\mathcal{L}_2}}{\|W(t)\|_{\mathcal{L}_2}}$ for this scenario in time-domain, listed as $\gamma_1 = 0.0226, \gamma_2 = 0.0234, \gamma_3 = 0.0166, \gamma_4 = 0.0187$. These results clearly agree with the desired performance  $\gamma_d = 1$, which validates the effectiveness of Theorem \ref{theo:2} for this scenario.

\subsection{Simulations with realistic vehicle dynamics}

In practice, each vehicle in a platoon is usually controlled by using either the throttle or the braking system according to a hierarchical architecture, consisting of an upper-level and a lower-level controller~\cite{li2015overview,rajamani2011vehicle}. The upper-level controller determines the desired acceleration using the information of its neighbors, which is the main topic of this paper, while the lower-level one generates the throttle or brake commands to track the desired acceleration trajectory.  The theoretical results of this paper work for the upper-level control stability, and it is better to test robustness using a realistic model of nonlinear vehicle dynamics. In this section, we present simulations of the platoon behavior in the presence of nonlinear vehicle dynamics that have not been explicitly considered in the process of the upper-level controller design.

As used in~\cite{wang2015longitudinal, Zheng2016distributed,dunbar2012distributed}, we consider the following model to describe longitudinal vehicle dynamics
\begin{equation} \label{eq:NonlinearVehicle}
    \begin{aligned}
        &\dot{p}_i = v_i(t), \dot{v}_i(t) = a_i(t), \\
        &a_i(t) = \frac{1}{m_i} \left(\eta_{T,i}\frac{T_i} {r_i} - C_{A,i} v_i^2 - m_igf \right), \\
        &\tau_i\dot{T}_i(t) + T_i(t) = T_{des,i}(t),
    \end{aligned}
\end{equation}
where $m_i$ is the vehicle mass, $C_{A,i}$ is the lumped aerodynamic drag coefficient, $g$ is the acceleration of gravity, $f$ is the coefficient of rolling resistance, $T_i(t)$ denotes the actual driving/braking torque, $T_{des,i}(t)$ is the desired driving/braking torque, $\tau_i$ is the inertial delay of vehicle longitudinal dynamics, $r_i$ denotes the tire radius, and $\eta_{T,i}$ is the mechanical efficiency of the driveline. In~\eqref{eq:NonlinearVehicle}, it is assumed that the powertrain dynamics are lumped to be a first-order inertial transfer function, which is widely used in the literature. Based on the desired acceleration that is calculated by the upper-level controller~\eqref{eq:6}, we employ an inverse model to generate the desired driving/braking torque
\begin{equation} \label{eq:InverseModel}
    T_{des,i}(t) = \frac{1}{\eta_{T,i}}\left(m_ia_{des,i} + C_{A,i}v_i^2 + m_igf\right)r_i.
\end{equation}

\begin{table}[t]
    \centering
    \renewcommand\arraystretch{1.2}
    \caption{Vehicle parameters used in the simulations}
    \label{tab:Nonlinear}
    \begin{tabular}{ m{1.1cm}<{\centering}|  m{0.3cm}<{\centering} m{0.3cm}<{\centering} m{0.3cm}<{\centering} m{0.3cm}<{\centering} m{0.3cm}<{\centering} m{0.3cm}<{\centering} m{0.3cm}<{\centering} m{0.3cm}<{\centering} m{0.3cm}<{\centering} m{0.3cm}<{\centering}}
        \hline \toprule[1pt]
        Number   & 1 & 2&  3&  4& 5 & 6 & 7 & 8 & 9 & 10\\\hline
         \begin{tabular}[x]{@{}c@{}}$m_i$\\[-0.25em]($\times 10^3$ kg)\end{tabular}   & 2.81 & 2.90 &  2.12&  2.91& 2.63 & 2.09 & 2.27 & 2.54 & 2.95 & 2.96\\
        $\tau_i$ (s) & 0.58 & 0.59&  0.51&  0.59& 0.56 & 0.50 & 0.52 & 0.55 & 0.60 & 0.60 \\
       % \hline
        \bottomrule[1pt]
        \end{tabular}
\end{table}

\begin{figure}
    \centering

    \subfigure[ ]
    {   \includegraphics[width=0.44\columnwidth]{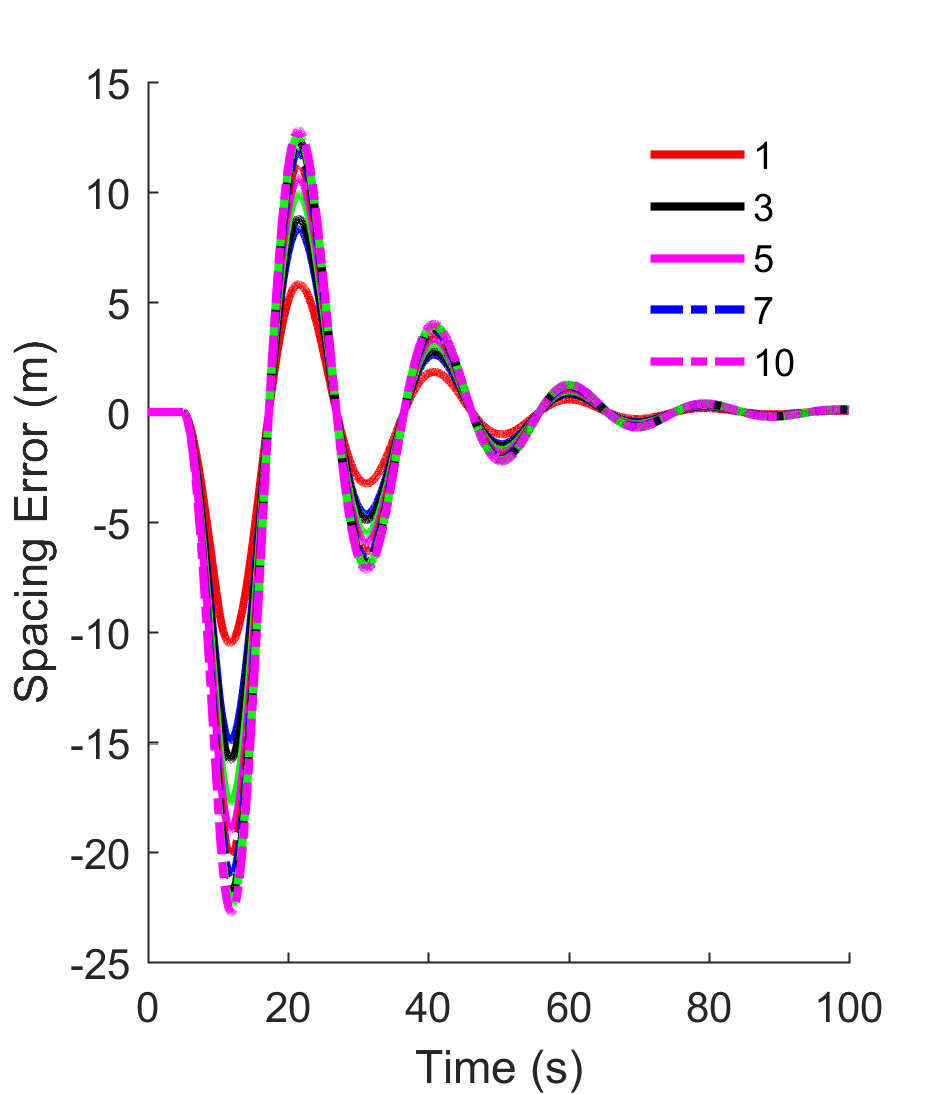}
    }
    \subfigure[ ]
    {
      \includegraphics[width=0.44\columnwidth]{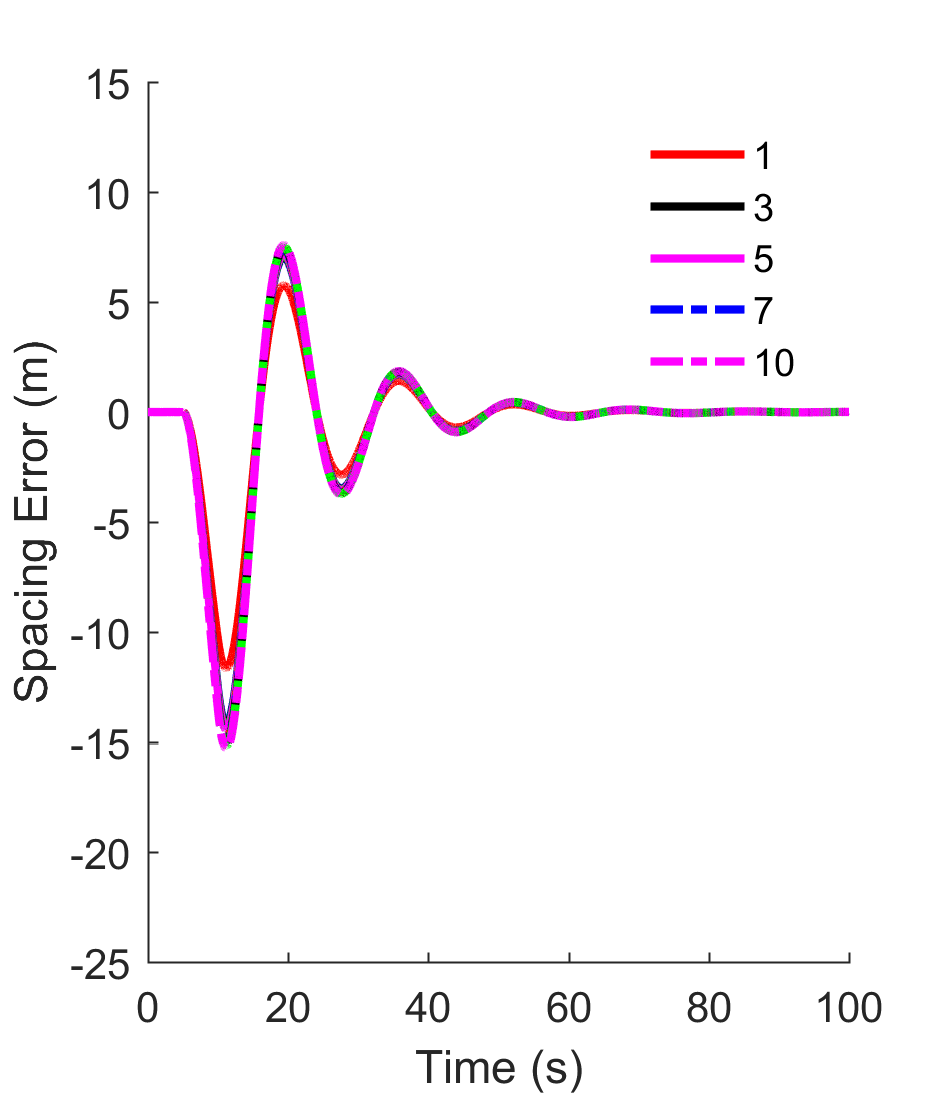}
    }
    \vspace{-4mm}
     \subfigure[ ]
    {
        \includegraphics[width=0.44\columnwidth]{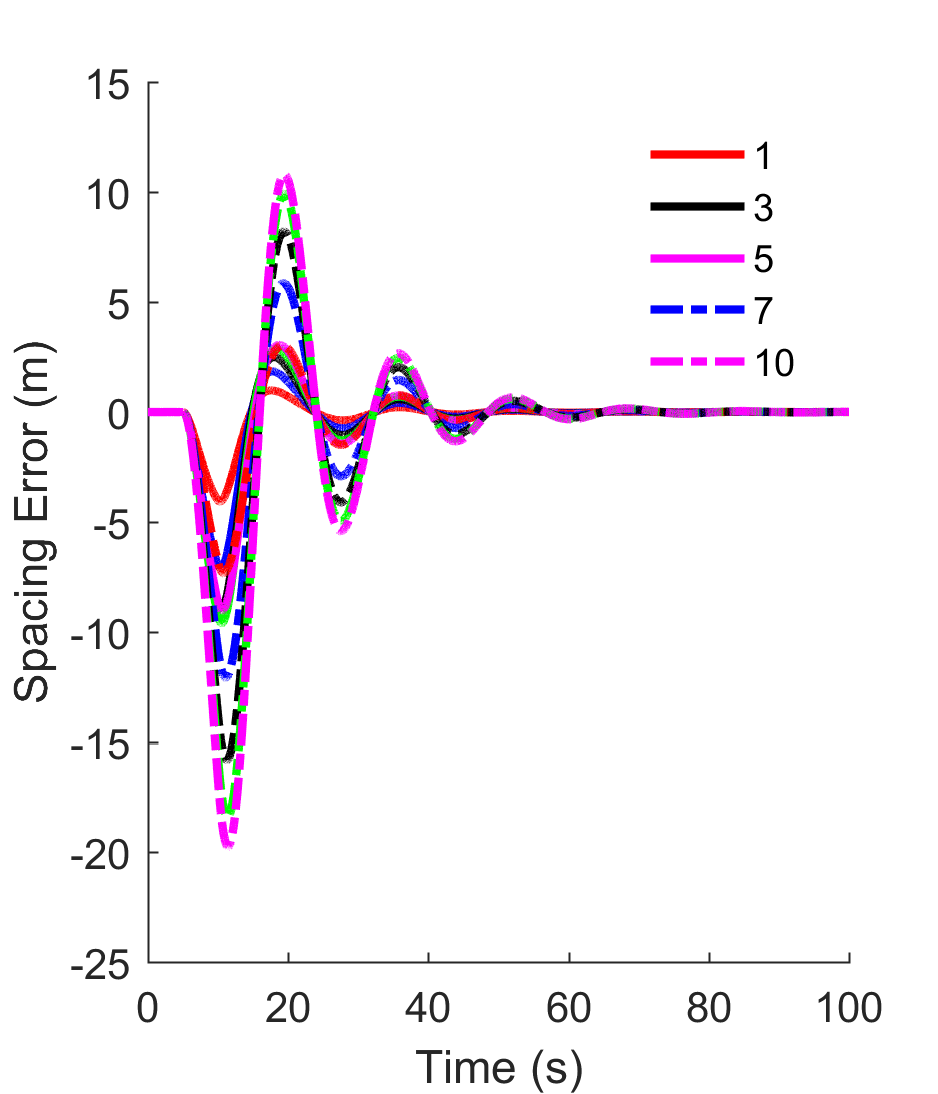}
    }
    \subfigure[ ]
    {
      \includegraphics[width=0.44\columnwidth]{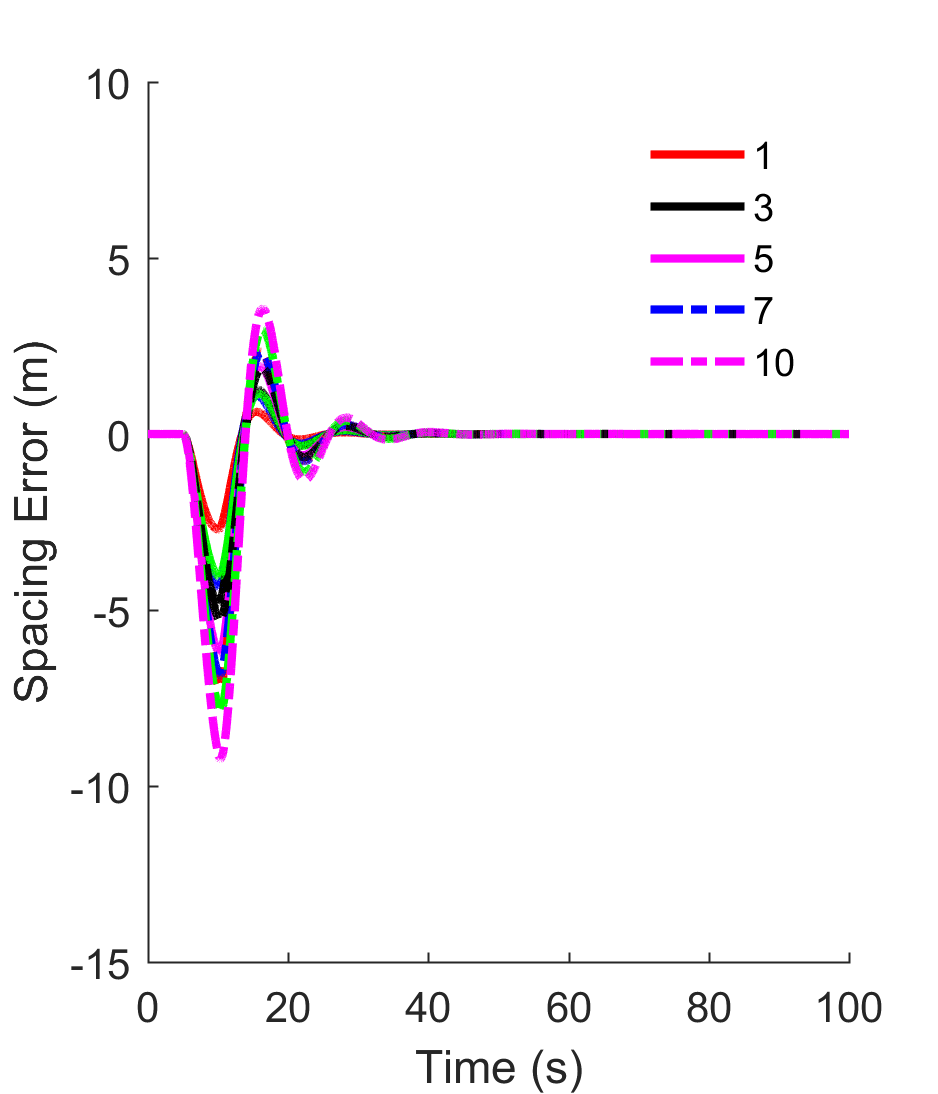}
    }

    \caption{Spacing errors for platoons ($N=10$) with the nonlinear vehicle dynamics~\eqref{eq:NonlinearVehicle} and the controller~\eqref{eq:NumeController}: (a) $h$-neighbor undirected topology with $h = 2$; (b) $h$-neighbor undirected topology with $h = 4$; (c) two multiple mini-platoons (both sizes are 5); (d) three multiple mini-platoons (the sizes are 3, 4 and 3).}
    \label{fig:NonlinearPlatoon}
\end{figure}

Then, we investigate the platoon behavior under the distributed controller~\eqref{eq:6} in the presence of nonlinear vehicle dynamics~\eqref{eq:NonlinearVehicle} (\emph{i.e.}, road friction, approximated engine dynamics, aerodrag forces). Similar to Section~\ref{section:SimulationII}, we consider a platoon of 11 vehicles, including one leader and ten followers. The communication topologies are the same with those in Section~\ref{section:SimulationII}.  The parameters of each following vehicle are randomly selected according to the passenger vehicles~\cite{wang2015longitudinal}. The value in our experiments is listed in Table~\ref{tab:Nonlinear}, and other parameters are $\eta_{T,i} = 0.9$, $g = 0.98 \; m/s^2$, $f_i = 0.01$, and $C_{A,i} = 0.492$. The acceleration or deceleration of the leader can be viewed as disturbances in a platoon. Motivated by~\cite{zheng2016stability}, we consider the following scenario: the initial state of the leader is set to $p_0(t) = 0, v_0 =
20\; m/s$, and the desired trajectory is given by
$$
    v_0(t) = \begin{cases} 20 \;m/s \qquad \quad t < 5\;s \\ 20+2t \;m/s \quad  5\;s \leq t \leq 10\;s \\ 30\; m/s \qquad \quad t> 10\;s \end{cases}.
$$

In the simulations, the desired spacing is set to $d_{i-1,i} = 25 \; m$, and the initial spacing errors and velocity errors are all equal to 0. The controller gain~\eqref{eq:NumeController} was used. Fig.~\ref{fig:NonlinearPlatoon} shows the spacing errors of the nonlinear platoon with different communication topologies. It is easy to see that the distributed controller~\eqref{eq:6} can stabilize the platoon with nonlinear vehicle dynamics. Furthermore, based on the theoretical analysis in Section~\ref{section:communication}, we know that the minimum eigenvalue of $\mathcal{L}+\mathcal{P}$ plays an important role on the performance of linear platoons. This statement also agrees with the results in Fig.~\ref{fig:NonlinearPlatoon} where nonlinear vehicle dynamics were considered. In fact, since the topology d) in Table~\ref{tab:1} has the largest $\lambda_{min}$, the transient performance in Fig.~\ref{fig:NonlinearPlatoon} (d) is fastest and the maximum spacing error is the smallest. These facts indicate that the theoretical analysis for linear platoons might be promising to serve as a guideline for the design of heterogeneous nonlinear platoons in practice.

\balance
\section{Conclusion}\label{section:conclusion}

This paper studies the robustness analysis and distributed $\mathcal{H}_{\infty}$ controller synthesis for platooning of connected vehicles with undirected topologies. Unified models in both time and frequency domain are derived to describe the collective behavior of homogeneous platoons with external disturbance using graph theory.
The major strategy of this paper is to equivalently decouple the collective dynamics of a platoon into a set of subsystems by exploiting the spectral decomposition of $\mathcal{L}+\mathcal{P}$. Along with this idea, we have derived the decomposition of platoon dynamics in both frequency domain (see \eqref{eq:29}) and time domain (see \eqref{eq:35}). Therefore, the robustness measure $\gamma$-gain is explicitly analyzed (see Theorem \ref{theo:1}), and the distributed $\mathcal{H}_{\infty} $ controller is also easily synthesized (see Theorem \ref{theo:2}). These results not only offer some insightful understandings of certain performance limits for large-scale platoons, but also hint us that coordination of multiple mini-platoons is one reasonable architecture to control large-scale platoons.

 One future direction is to address the influence of imperfect communication, such as time delay and packet loss, on the robustness performance. Besides, this paper exclusively focus on the identical feedback gains \eqref{eq:6}. Whether or not non-identical controllers can improve the robustness performance of a platoon is an interesting topic for future research. We notice that some recent work proposed asymmetric controller for platoons with BD topology, resulting in certain improvements in terms of stability margin (see, \emph{e.g.}, \cite{barooah2009mistuning, zheng2015stabilityMargin,tangerman2012asymmetric}).

\bibliographystyle{IEEEtran}
\bibliography{Reference}

\end{document}